\documentclass[3p]{elsarticle}

\usepackage{lmodern}
\usepackage{amsmath}
\usepackage{amssymb}
\usepackage{float}
\usepackage{mathrsfs}
\usepackage{tikz}
\usepackage{tikz-cd}
\usepackage{graphics}
\usepackage{booktabs}
\usepackage{multirow}
\usepackage{hyperref}
\usepackage{siunitx}
\usepackage[font=small,labelfont=bf]{caption}
\usepackage [english]{babel}
\usepackage{url}
\usepackage{array}
\usepackage{booktabs}

\newtheorem{theorem}{Theorem}
\newtheorem{proposition}{Proposition}
\newtheorem{corollary}{Corollary}
\newtheorem{lemma}{Lemma}
\newtheorem{mainthm}{Theorem}

\newdefinition{example}{Example}
\newdefinition{definition}{Definition}
\newdefinition{remark}{Remark}
\newproof{proof}{Proof}

\newcommand{\Ric}{\operatorname{Ric}}
\newcommand{\Hess}{\operatorname{Hess}}
\newcommand{\dv}{\operatorname{div}}
\newcommand{\grad}{\operatorname{grad}}
\newcommand{\scal}{\operatorname{scal}}
\newcommand{\tr}{\operatorname{tr}}
\newcommand{\Span}{\operatorname{span}}
\newcommand{\Id}{\mathrm{Id}}
\newcommand{\R}{\mathbb{R}}
\newcommand{\Lie}{\mathcal{L}}
\newcommand{\eps}{\varepsilon}
\newcommand{\sym}{\odot}
\newcommand{\nt}{\tilde\nabla}
\begin{document}
    \title{$C_{12}$- and para-$C_{12}$-manifolds: twisted products, $\eta$-Einstein structures and Ricci solitons of the $\eta$-reversed metric}

    \author[a]{Umut Selvi}
        \author[a]{Derya Sağlam\corref{cor1}}
            \ead{derya.saglam@hbv.edu.tr}
            \cortext[cor1]{Corresponding author}

    \address[a]{Department of Mathematics, Ankara Hacı Bayram Veli University, Polatlı Faculty of Science and Art, 06900, Ankara, Turkey}

    \begin{abstract}
        The class $C_{12}$ of Chinea and Gonz\'alez consists of almost contact metric manifolds which are locally twisted products $h+\sigma^2dt^2$ over K\"ahler manifolds. We show that the paracontact analogue $\mathbb G_{12}$ consists, in the same way, of twisted products over para-K\"ahler manifolds, and we relate $g$ to the $\eta$-reversed metric $\tilde g=g-2\eta\otimes\eta$, which is Lorentzian in the almost contact case. For $\eta$-Einstein structures with constant $\lambda$ the twisting function satisfies, on the base, the quasi-Einstein equation with $m=1$; using the structure theory of quasi-Einstein metrics on K\"ahler manifolds of arbitrary signature and on para-K\"ahler manifolds, we prove that either $g$ is Einstein or $\omega=-\nabla_\xi\eta$ is exact, give a local classification (complete in the almost contact case, partial in the paracontact case), and show that compact $\eta$-Einstein $C_{12}$-manifolds with constant $\lambda$ are normal. Ricci solitons $(\tilde g,f\xi,\lambda)$ either have a Killing potential or are locally products of an Einstein manifold with a timelike line. We also classify the left-invariant $C_{12}$-structures on three-dimensional Lie groups and correct several formulas in the literature.
    \end{abstract}

    \begin{keyword}
        $C_{12}$-manifold; paracontact metric manifold; twisted product; $\eta$-Einstein manifold; Ricci soliton; Lorentzian metric; quasi-Einstein metric
        \MSC[2020] 53C25; 53C50; 53C15; 53D15; 53E20
    \end{keyword}

    \maketitle

    \section{Introduction}
    \noindent
    The class $C_{12}$ of Chinea and Gonz\'alez \cite{CG90} is one of the twelve basic classes of almost contact metric manifolds; its paracontact analogue is the class $\mathbb{G}_{12}$ of \cite{ZN18}, whose members we call para-$C_{12}$-manifolds. A $C_{12}$-manifold carries a Reeb field $\xi$ with $\nabla_X\xi=-\eta(X)\psi$, where $\psi=-\nabla_\xi\xi$, and the $1$-form $\omega=g(\psi,\cdot)$ satisfies $d\eta=\omega\wedge\eta$. The members of $C_{12}$ are always integrable but normal only in the trivial case $\omega=0$, and they are locally twisted products $g=h+\sigma^2dt^2$ over K\"ahler manifolds \cite{dCF19,Fal13}; three-dimensional $C_{12}$-manifolds have recently been studied from several viewpoints \cite{BB22,Bel24,BBB21,dCF19}. We show that para-$C_{12}$-manifolds are, in the same way, twisted products over para-K\"ahler manifolds (Theorem~\ref{thm:twisted}). Reversing the sign of $g$ on $\xi$ gives the \emph{$\eta$-reversed metric} $\tilde g=g-2\eta\otimes\eta$, the member $t=-2$ of the canonical variation of Olea \cite{Ole}; it is Lorentzian in the almost contact case and has signature $(n,n+1)$ in the paracontact case. Delloum and Beldjilali \cite{DB24} initiated the study of Ricci solitons for the pair $(g,\tilde g)$; several of their formulas are not correct, and we list the corrections in Section~\ref{sec:relation}. We relate the geometries of $g$ and $\tilde g$ without using the twisted model (Theorem~\ref{thm:compare}). For an $\eta$-Einstein structure the twisting function satisfies, on the base, the quasi-Einstein equation $\Ric_h-\frac1\sigma\Hess_h\sigma=\lambda h$, i.e.\ the case $m=1$ of \eqref{eq:QE}; we use the structure theory of quasi-Einstein metrics on K\"ahler manifolds of arbitrary signature and on para-K\"ahler manifolds developed in the companion paper \cite{SS-QE}, and we recall what we need in Section~\ref{sec:QErecall}. Throughout, \emph{$\eta$-Einstein with constant $\lambda$} means that $S(X,Y)=\lambda g(X,Y)$ for all $X,Y\perp\xi$ with a constant $\lambda$; the coefficient of $\eta\otimes\eta$ is then in general not constant (Definition~\ref{def:etaE}), and without the constancy of $\lambda$ the results below fail (Remark~\ref{rem:nonconst}).

    \begin{mainthm}\label{thmB}
        Let $(M^{2n+1},\varphi,\xi,\eta,g)$ be a connected $C_{12}$-manifold which is $\eta$-Einstein with constant $\lambda$.
        \begin{itemize}
            \item[(1)] Either $g$ is Einstein, $S=\lambda g$, or $\omega$ is exact on $M$. In particular, if $d\omega\neq0$ at one point, then $g$ is Einstein.
            \item[(2)] Every point of a dense open subset of $M$ has a neighbourhood which, as a $C_{12}$-manifold, is of one of the following types: (a) normal, $g=h+dt^2$ with $\Ric_h=\lambda h$; (b) the product of a three-dimensional $\eta$-Einstein $C_{12}$-manifold $(\Sigma\times I,h_\Sigma+a^2dt^2)$, $a=a(\Sigma)$, $da\neq0$, and a K\"ahler--Einstein manifold $(K,h_K)$ with $\Ric_K=\lambda h_K$; (c) $\lambda\neq0$ and the product of a three-dimensional space form of curvature $\lambda/2$, carrying a $C_{12}$-structure with non-closed $\omega$, and a K\"ahler--Einstein manifold with $\Ric_K=\lambda h_K$; (d) $\lambda=0$ and the product of a flat $C_{12}$-manifold with non-closed $\omega$ and a Ricci-flat K\"ahler manifold. In (b)--(d) the K\"ahler factor is absent if $n=1$.
            \item[(3)] If $M$ is compact, then the structure is normal: $\omega=0$, $\xi$ is parallel, $S=\lambda(g-\eta\otimes\eta)$, and locally $g=h+dt^2$ with $\Ric_h=\lambda h$. In particular, $g$ is Einstein only if it is Ricci-flat.
        \end{itemize}
    \end{mainthm}

    \begin{mainthm}\label{thmBP}
        Let $(M^{2n+1},\varphi,\xi,\eta,g)$ be a connected para-$C_{12}$-manifold which is $\eta$-Einstein with constant $\lambda$. Then part (1) of Theorem~\ref{thmB} holds, and if $M$ is compact, then $\lambda=0$ implies $S=0$, while $\lambda\neq0$ implies that $g$ is not Einstein and that $\omega$ is exact. Where $\omega$ is closed and $\psi$ is not null, $M$ splits locally as in Theorem~\ref{thmB}(2)(b), with $\Sigma$ a Lorentzian surface and $K$ para-K\"ahler--Einstein. Where $\omega$ is closed and $\psi$ is null and non-zero, $g$ and $\tilde g$ are Ricci-flat and, if both components of $\psi$ in the eigendistributions of $\varphi$ are non-zero, no three-dimensional factor splits off. Where $d\omega\neq0$, either the analogue of type (c) holds, or $\lambda=0$ and $g$ is Ricci-flat; in the latter case we obtain no normal form.
    \end{mainthm}

    \noindent
    Theorems~\ref{thmB} and~\ref{thmBP} are proved in Section~\ref{sec:etaE}: part (1) is Theorem~\ref{thm:nonclosed}, part (2) of Theorem~\ref{thmB} is Corollary~\ref{cor:classR}, which rests on Theorem~\ref{thm:nonclosedR} and Corollary~\ref{cor:R}, the compact case is Theorem~\ref{thm:compact}, and the paracontact statements are Theorem~\ref{thm:nonclosedR}(b) and Corollary~\ref{cor:P}. The normality in the compact case (Theorem~\ref{thm:compact}(c)) is obtained by combining the first integral of Lemma~\ref{lem:KK} with a Myers-type argument on the leaves of $\ker\eta$ and with the Killing field $P\nabla u$ of Theorem~\ref{thm:QE}; it uses the Riemannian signature, and the corresponding paracontact question is open (Section~\ref{sec:open}). For Ricci solitons whose potential is collinear with $\xi$ we obtain a global result; it holds for every unit vector field $\xi$ with $\nabla\xi=-\eta\otimes\psi$, i.e.\ for every twisted product with a one-dimensional fibre (setting (Tw) of Remark~\ref{rem:settingE}).

    \begin{mainthm}\label{thmA}
        Let $(M,\varphi,\xi,\eta,g)$ be a connected $C_{12}$- or para-$C_{12}$-manifold and let $(\tilde g,f\xi,\lambda)$ be a Ricci soliton. Then either $\lambda\neq0$, $\xi$ is parallel, and locally $\tilde g=-d\tau^2+h$ with $\Ric_h=\lambda h$ and $f\xi=(\lambda\tau+k_0)\partial_\tau$; or $g$ and $\tilde g$ are Einstein with constant $\lambda$ and $f\xi$ is a Killing field of $\tilde g$. In the second case $\omega$ is exact unless $f\equiv0$.
    \end{mainthm}

    \noindent
    For potentials $V$ orthogonal to $\xi$, the Ricci solitons $(\tilde g,V,\lambda)$ and $(g,V,\lambda)$ coincide (Corollary~\ref{cor:Vperp}); this gives non-trivial solitons of $\tilde g$, for instance from Hamilton's cigar (Example~\ref{ex:cigar}).

    \paragraph{Relation to the literature} Formula (i) of Theorem~\ref{thm:compare} is a special case of the expression for the Levi-Civita connection of $g+t\,\eta\otimes\eta$ in \cite[Prop.~2.3]{Ole}, and the remaining identities of Theorem~\ref{thm:compare} follow from it by elementary computations (for warped products, (ii) also follows from \cite[Cor.~7.43]{ONe83}); the first integral of Lemma~\ref{lem:KK} is due to Kim and Kim \cite{KK03}, and its proof does not use the signature. Theorem~\ref{thm:nonclosed}(i) is, in arbitrary signature, the case $m=1$ of the analysis of the quadratic form $\mu$ on the space of solutions of the quasi-Einstein equation by He, Petersen and Wylie \cite{HPW15b}; equivalently, where $\omega$ is not closed the functions $\sigma(\cdot,t)$ are static potentials in the sense of Fischer and Marsden \cite{FM75} and Corvino \cite{Cor00} (Remark~\ref{rem:HPW}). The new points in Theorems~\ref{thmB} and~\ref{thmBP} are the global identity $d\Theta=-2\Theta\omega$, the use of the (para-)K\"ahler structure of the base through the results of \cite{SS-QE}, which forces a two-dimensional factor, the isotropic paracontact case, and the normality of compact $\eta$-Einstein $C_{12}$-manifolds with constant $\lambda$, whose proof uses, for $\lambda>0$, the case $m=1$ of a Myers-type theorem for the Bakry--\'Emery Ricci tensor \cite{Qian97}, in a form which we derive from the classical Myers theorem on a warped product. Ricci solitons on warped products have mainly been studied for gradient potentials \cite{BT22,FFG17}; the potentials $f\xi$ in Theorem~\ref{thmA} are in general not gradients, and Theorem~\ref{thmA} shows that the twisting function must be a warping function. To the best of our knowledge, the twisted-product description of para-$C_{12}$-manifolds and Theorems~\ref{thmB}, \ref{thmBP} and~\ref{thmA} do not appear in the literature.

    \paragraph{Organisation} Section~\ref{sec:prelim} fixes conventions and recalls the results of \cite{SS-QE} on quasi-Einstein metrics that we use. Section~\ref{sec:twisted} contains the twisted-product description of $C_{12}$- and para-$C_{12}$-manifolds, and Section~\ref{sec:compare} the comparison of $g$ and $\tilde g$. Sections~\ref{sec:etaE} and~\ref{sec:solitons} treat $\eta$-Einstein structures and Ricci solitons. Examples are collected in Section~\ref{sec:examples}, and Section~\ref{sec:lie} treats left-invariant structures on three-dimensional Lie groups. Section~\ref{sec:relation} summarises the corrections to \cite{DB24}, whose details are given in Appendix~\ref{app:DB}, and Section~\ref{sec:open} lists open questions.
    \section{Preliminaries}\label{sec:prelim}
    \noindent
    For $1$-forms we use $d\eta(X,Y)=X\eta(Y)-Y\eta(X)-\eta([X,Y])$ and $(\alpha\wedge\beta)(X,Y)=\alpha(X)\beta(Y)-\alpha(Y)\beta(X)$. The curvature tensor is $R(X,Y)=[\nabla_X,\nabla_Y]-\nabla_{[X,Y]}$ and the Ricci tensor of a pseudo-Riemannian metric $g$ is $S(Y,Z)=\tr\bigl(X\mapsto R(X,Y)Z\bigr)$; $r$ denotes the scalar curvature. For other metrics we use the corresponding notation, e.g.\ $\tilde R$, $\tilde S$, $\tilde r$ for $\tilde g$, and $\Ric_h$, $\Hess_h$, $\Delta_h=\tr_h\Hess_h$ and $\scal_h$ for a metric $h$. For $1$-forms $\alpha,\beta$ we put $\alpha\sym\beta:=\alpha\otimes\beta+\beta\otimes\alpha$; for a vector $X$ we write $X^\flat$ for the metrically dual $1$-form, $\alpha^\sharp$ for the vector dual to a $1$-form $\alpha$, and $(v\otimes\alpha)(X):=\alpha(X)v$. All manifolds are connected, and all results are local unless stated otherwise.

    \subsection{Quasi-Einstein metrics on K\"ahler and para-K\"ahler manifolds}\label{sec:QErecall}
    \noindent
    Let $(N,h)$ be a pseudo-Riemannian manifold, $u>0$ a smooth function, $m\in\R\setminus\{0\}$ and $\lambda\in\R$. The \emph{quasi-Einstein equation} is
    \begin{equation}\label{eq:QE}
        \Ric_h-\frac{m}{u}\Hess_hu=\lambda h .
    \end{equation}
    For $m=1$ it is satisfied by the twisting function of an $\eta$-Einstein $C_{12}$- or para-$C_{12}$-manifold with constant $\lambda$ (Section~\ref{sec:etaE}). We recall the results of \cite{SS-QE} that we need.

    Kim and Kim \cite{KK03} state this for Riemannian metrics; since we need it in arbitrary signature, we include the short proof.
    
    \begin{lemma}\label{lem:KK}
        Let $u>0$ satisfy \eqref{eq:QE}. Then $\mu:=u\,\Delta_hu+(m-1)\,h(\nabla u,\nabla u)+\lambda u^2$ is locally constant.
    \end{lemma}

    \noindent
    The proof in \cite{KK03} (see also \cite[(2.6)]{CSW11}) is a pointwise computation based on the identities $\dv\Ric_h=\frac12d\scal_h$ and $\dv\Hess_hu=d\Delta_hu+\Ric_h(\nabla u,\cdot)$, which hold in every signature; hence the lemma holds for pseudo-Riemannian metrics of arbitrary signature \cite[Lemma~2.1]{SS-QE}.

    Let $\eps\in\{-1,1\}$. An \emph{$\eps$-K\"ahler manifold} is a pseudo-Riemannian manifold $(N^{2n},h)$ with an endomorphism field $P$ such that $P^2=\eps\,\Id$, $h(PX,PY)=-\eps\,h(X,Y)$ and $\nabla P=0$. For $\eps=-1$ this is a K\"ahler manifold of arbitrary signature; for $\eps=1$ it is a para-K\"ahler manifold \cite{AMT09,CFG96}, $h$ has signature $(n,n)$, and the eigendistributions $T^\pm:=\ker(P\mp\Id)$ are parallel and totally null.

    \begin{lemma}[{\cite[Lemma~3.1(a)]{SS-QE}}]\label{lem:Ricform}
        On an $\eps$-K\"ahler manifold, $\Ric_h(PX,PY)=-\eps\Ric_h(X,Y)$.
    \end{lemma}

    \begin{theorem}[{\cite[Theorem~3.2]{SS-QE}}]\label{thm:QE}
        Let $(N^{2n},h,P)$ be $\eps$-K\"ahler, let $u>0$ satisfy \eqref{eq:QE}, and put $\zeta:=\nabla u$ and $A:=\nabla\zeta$, so that $h(AX,Y)=\Hess_hu(X,Y)$.
        \begin{itemize}
            \item[(i)] $AP=PA$, and $P\zeta$ is a Killing vector field.
            \item[(ii)] On an open set where $h(\zeta,\zeta)\neq0$, the distribution $\mathcal Z:=\Span\{\zeta,P\zeta\}$ has rank two, is non-degenerate and parallel, and $A(\mathcal Z)\subset\mathcal Z$, $A(\mathcal Z^\perp)=0$. Consequently every point of this set has a neighbourhood isometric, as an $\eps$-K\"ahler manifold, to a product $(\Sigma^2,h_\Sigma,P_\Sigma)\times(K^{2n-2},h_K,P_K)$ with $T\Sigma=\mathcal Z$, $u=u(\Sigma)$, $\Ric_K=\lambda h_K$ and $\Ric_\Sigma-\frac mu\Hess_\Sigma u=\lambda h_\Sigma$.
            \item[(iii)] Let $\eps=1$, let $U$ be a connected open set on which $\zeta\neq0$ and $h(\zeta,\zeta)=0$, and write $\zeta=\zeta_++\zeta_-$ with $\zeta_\pm\in T^\pm$. Then $\lambda=0$ and $\Delta_hu=0$ on $U$, and every point of a dense open subset of $U$ has a neighbourhood on which one of the following holds:
            \begin{itemize}
                \item[(iii-1)] $\zeta\in T^+$ or $\zeta\in T^-$; then $\Hess_hu=0$, $\zeta$ is a parallel null vector field and $\Ric_h=0$;
                \item[(iii-2)] $\zeta_+\neq0\neq\zeta_-$; then $\Hess_hu=c\,\zeta_+^\flat\sym\zeta_-^\flat$ for a smooth function $c$, and every non-degenerate $P$-invariant subspace of $T_pN$ containing $\zeta_p$ has dimension at least four.
            \end{itemize}
        \end{itemize}
    \end{theorem}

    \noindent
    For Riemannian K\"ahler metrics, (ii) is the local form of \cite[Thm.~1.3]{CSW11}. Normal forms for the solutions of type (iii-2), and for the analogous isotropic solutions in the pseudo-K\"ahler case, are given in \cite[Theorems~4.1 and~5.1]{SS-QE}.

    \section{$C_{12}$- and para-$C_{12}$-manifolds}\label{sec:twisted}
    \noindent
    We work in one of two settings.
    \begin{itemize}
        \item[(R)] \emph{Riemannian:} $(\varphi,\xi,\eta,g)$ is an almost contact metric structure on $M^{2n+1}$, i.e.\ $\varphi^2=-\Id+\eta\otimes\xi$, $\eta(\xi)=1$, $g(\varphi X,\varphi Y)=g(X,Y)-\eta(X)\eta(Y)$ and $g(X,\xi)=\eta(X)$.
        \item[(P)] \emph{Paracontact:} $(\varphi,\xi,\eta,g)$ is an almost paracontact metric structure on $M^{2n+1}$, i.e.\ $\varphi^2=\Id-\eta\otimes\xi$, $\eta(\xi)=1$, $g(\varphi X,\varphi Y)=-g(X,Y)+\eta(X)\eta(Y)$ and $g(X,\xi)=\eta(X)$; the signature of $g$ is $(n+1,n)$, Lorentzian for $n=1$.
    \end{itemize}
    In both settings $\varphi\xi=0$, $\eta\circ\varphi=0$, $g(\varphi X,Y)=-g(X,\varphi Y)$, $g(\xi,\xi)=1$, and $D:=\ker\eta=\varphi(TM)$ carries a non-degenerate metric. Put $F(X,Y,Z)=g((\nabla_X\varphi)Y,Z)$ and $\Phi(X,Y)=g(X,\varphi Y)$.

    \begin{definition}\label{def:C}
        The structure satisfies \emph{condition (C)} if
        \begin{equation}\label{eq:C}
            F(X,Y,Z)=\eta(X)\bigl\{\eta(Y)F(\xi,\xi,Z)-\eta(Z)F(\xi,\xi,Y)\bigr\}.
        \end{equation}
        A \emph{$C_{12}$-manifold} (resp.\ \emph{para-$C_{12}$-manifold}) is a manifold with a structure satisfying (C) in setting (R) (resp.\ (P)). We put $\psi:=-\nabla_\xi\xi$, $\omega:=g(\psi,\cdot)$ and $|\psi|^2:=g(\psi,\psi)$; in setting (P), $|\psi|^2$ may have either sign.
    \end{definition}

    \noindent
    Since $(\nabla_X\Phi)(Y,Z)=F(X,Z,Y)$ and $F(\xi,\xi,Z)=g(\varphi\psi,Z)=(\nabla_\xi\eta)(\varphi Z)$, condition (C) is equivalent to $(\nabla_X\Phi)(Y,Z)=\eta(X)\eta(Z)(\nabla_\xi\eta)\varphi Y-\eta(X)\eta(Y)(\nabla_\xi\eta)\varphi Z$, which is the defining relation of the class $C_{12}$ in \cite{CG90}; in setting (P), (C) is the defining relation of the class $\mathbb{G}_{12}$ of \cite{ZN18}.

    \begin{lemma}\label{lem:basic}
        Under condition (C):
        \begin{itemize}
            \item[(i)] $\nabla_X\xi=-\eta(X)\psi$; hence $(\nabla_X\eta)Y=-\eta(X)\omega(Y)$, $\omega(\xi)=0$ and $\eta(\psi)=0$;
            \item[(ii)] $d\eta=\omega\wedge\eta$;
            \item[(iii)] $(\nabla_X\varphi)Y=\eta(X)\bigl\{\omega(\varphi Y)\xi+\eta(Y)\varphi\psi\bigr\}$;
            \item[(iv)] $S(X,\xi)=-\eta(X)\dv\psi$ and $\eta(\nabla_Y\psi)=\eta(Y)|\psi|^2$;
            \item[(v)] $N_\varphi=0$ and $d\Phi=0$, where $N_\varphi(X,Y)=\varphi^2[X,Y]+[\varphi X,\varphi Y]-\varphi[\varphi X,Y]-\varphi[X,\varphi Y]$.
        \end{itemize}
    \end{lemma}

    \begin{proof}
        (i) Putting $Y=\xi$ in \eqref{eq:C} gives $F(X,\xi,Z)=\eta(X)F(\xi,\xi,Z)$. Since $(\nabla_X\varphi)\xi=-\varphi\nabla_X\xi$ and $\varphi$ is $g$-skew, $F(X,\xi,Z)=g(\nabla_X\xi,\varphi Z)$. Hence $W:=\nabla_X\xi-\eta(X)\nabla_\xi\xi$ satisfies $g(W,\varphi Z)=0$ for all $Z$. As $g(\xi,\xi)=1$, $W\in D=\varphi(TM)$, and $g|_D$ is non-degenerate, so $W=0$.
        (ii) $d\eta(X,Y)=g(\nabla_X\xi,Y)-g(\nabla_Y\xi,X)=-\eta(X)\omega(Y)+\eta(Y)\omega(X)$.
        (iii) $F(\xi,\xi,Z)=-g(\varphi\nabla_\xi\xi,Z)=g(\varphi\psi,Z)=-\omega(\varphi Z)$; insert this into \eqref{eq:C}.
        (iv) From (i), $R(X,Y)\xi=\eta(X)\omega(Y)\psi-\eta(Y)\omega(X)\psi-\eta(Y)\nabla_X\psi+\eta(X)\nabla_Y\psi$. Taking the trace in $X$ and using $\eta(\psi)=0$, $\omega(\psi)=|\psi|^2$ and $\eta(\nabla_Y\psi)=-g(\psi,\nabla_Y\xi)=\eta(Y)|\psi|^2$ gives $S(Y,\xi)=-\eta(Y)\dv\psi$.
        (v) By (i) and (iii), $\nabla_X\varphi=\eta(X)\nabla_\xi\varphi$ and $\varphi(\nabla_X\varphi)Y=\eta(X)\eta(Y)\varphi^2\psi$. Since $\nabla$ is torsion-free,
        \[
            N_\varphi(X,Y)=(\nabla_{\varphi X}\varphi)Y-(\nabla_{\varphi Y}\varphi)X+\varphi(\nabla_Y\varphi)X-\varphi(\nabla_X\varphi)Y ;
        \]
        the first two terms vanish because $\eta\circ\varphi=0$, and the last two cancel. Finally, $(\nabla_X\Phi)(Y,Z)=g(Y,(\nabla_X\varphi)Z)=-F(X,Y,Z)$, and the cyclic sum of the right-hand side of \eqref{eq:C} over $X,Y,Z$ vanishes; hence $d\Phi=0$.
    \end{proof}

    \noindent
    In setting (P), $N_\varphi$ in (v) is given by the same expression. With our convention for $d$, an almost contact metric structure is normal if and only if $N_\varphi+d\eta\otimes\xi=0$, and an almost paracontact metric structure is normal if and only if $N_\varphi-d\eta\otimes\xi=0$ (in the convention $d\eta(X,Y)=\frac12(X\eta(Y)-Y\eta(X)-\eta([X,Y]))$ these conditions read $N_\varphi+2\,d\eta\otimes\xi=0$ and $N_\varphi-2\,d\eta\otimes\xi=0$, respectively). Thus every $C_{12}$- or para-$C_{12}$-manifold has $N_\varphi=0$ and $d\Phi=0$; by Theorem~\ref{thm:twisted} below it is normal only if $\omega=0$. Note also that $S(\xi,\xi)=-\dv\psi$. The identities of Lemma~\ref{lem:basic} are equivalent to those in \cite{BBB21,dCF19}, where they are stated in dimension three or obtained by other methods.

    \begin{definition}\label{def:twisted}
        Let $(N,h)$ be a pseudo-Riemannian manifold, $I\subset\R$ an open interval and $\sigma$ a positive function on $N\times I$. The \emph{twisted product} $N\times_\sigma I$ is $N\times I$ with the metric $h+\sigma^2dt^2$ \cite{PR93}.
    \end{definition}

    \begin{theorem}\label{thm:twisted}
        Let $(M^{2n+1},\varphi,\xi,\eta,g)$ be a $C_{12}$-manifold (setting (R)) or a para-$C_{12}$-manifold (setting (P)). Every point has a neighbourhood $N\times I$ with coordinates $(x,t)$ in which
        \begin{equation}\label{eq:model}
            g=h+\sigma(x,t)^2dt^2,\qquad \eta=\sigma\,dt,\qquad \xi=\sigma^{-1}\partial_t ,
        \end{equation}
        $h$ is a metric on $N$, $\varphi|_D$ is independent of $t$, and $(N,P:=\varphi|_D,h)$ is K\"ahler in (R) and para-K\"ahler in (P). The leaves $N\times\{t\}$ are the leaves of $D$ and the curves $\{x\}\times I$ are the integral curves of $\xi$; hence $(N,h,P)$ is determined by the structure up to isomorphism, and $\sigma$ is determined up to the reparametrisations $s=T(t)$, $T'>0$, under which $\sigma$ is replaced by $\sigma/T'(t)$. Conversely, for any such data the formulas \eqref{eq:model} together with $\varphi\partial_t=0$ and $\varphi|_{TN}=P$ define a $C_{12}$- (resp.\ para-$C_{12}$-) structure, with $\psi=\grad_h\ln\sigma$ and $\omega=d_N\ln\sigma$. Moreover $d\omega=0$ if and only if $\partial_t\partial_{x^i}\ln\sigma=0$ for all $i$, if and only if $\sigma=\alpha(t)a(x)$; after reparametrising $t$ one then has $\sigma=a(x)$. The structure is normal if and only if $\omega=0$.
    \end{theorem}

    \noindent
    In setting (R) the first assertion is due to Falcitelli \cite{Fal13} and de Candia and Falcitelli \cite{dCF19}; we give a proof that works in both settings.

    \begin{proof}
        \emph{Splitting.} By Lemma~\ref{lem:basic}(ii), $\eta\wedge d\eta=0$, so $D$ is integrable. For $X,Y\in D$ we have $g(\nabla_XY,\xi)=-g(Y,\nabla_X\xi)=0$, so the leaves of $D$ are totally geodesic. The foliation by integral curves of $\xi$ is one-dimensional, and every one-dimensional foliation is totally umbilic. Since $g(\xi,\xi)=1$, both distributions are non-degenerate, and they are orthogonal. By the theorem of Ponge and Reckziegel \cite{PR93}, which holds in every signature, $M$ is locally a twisted product $N\times_\sigma I$, that is, $g=h+\sigma^2dt^2$ with $h$ independent of $t$, where the factors are the leaves of the two foliations. Then $\eta=g(\xi,\cdot)=\sigma\,dt$ and $\xi=\sigma^{-1}\partial_t$ (after possibly replacing $t$ by $-t$). A coordinate $t$ with $\eta\in\R\,dt$ is constant on the leaves of $D$, so any other such coordinate with the same orientation is $s=T(t)$ with $T'>0$, and then $\sigma\,dt=(\sigma/T'(t))\,ds$.

        \emph{$\varphi|_D$ is independent of $t$.} For a function $u$ and $Y\in TM$, $(\Lie_{u\xi}\varphi)Y=u(\Lie_\xi\varphi)Y-(\varphi Y)(u)\,\xi$, because $\varphi\xi=0$. For $Y\in D$, $(\Lie_\xi\varphi)Y=(\nabla_\xi\varphi)Y-\nabla_{\varphi Y}\xi+\varphi\nabla_Y\xi=(\nabla_\xi\varphi)Y$ by Lemma~\ref{lem:basic}(i), and $(\nabla_\xi\varphi)Y=\omega(\varphi Y)\xi$ by Lemma~\ref{lem:basic}(iii). Hence $\Lie_{\partial_t}\varphi(\partial_{x^i})=\Lie_{\sigma\xi}\varphi(\partial_{x^i})\in\Span\{\xi\}$. On the other hand $\Lie_{\partial_t}\varphi(\partial_{x^i})=(\partial_t\varphi^j_i)\,\partial_{x^j}\in D$. Therefore $\partial_t\varphi|_D=0$.

        \emph{(Para-)K\"ahler base.} For $X,Y\in D$ we have $(\nabla_X\varphi)Y=0$ by Lemma~\ref{lem:basic}(iii), since $\eta(X)=0$. As the leaves of $D$ are totally geodesic, $\nabla_XY$ and $\nabla_X(\varphi Y)$ are the Levi-Civita derivatives along the leaf $N\times\{t\}$, which is isometric to $(N,h)$; hence $\nabla^hP=0$. Since $P^2=-\Id$ and $h(PX,PY)=h(X,Y)$ in (R), resp.\ $P^2=\Id$ and $h(PX,PY)=-h(X,Y)$ in (P), $(N,P,h)$ is K\"ahler, resp.\ para-K\"ahler.

        \emph{Converse.} Let $X,Y\in TN$ be $t$-independent lifts. Then $\nabla_XY=\nabla^h_XY$, $\nabla_X\partial_t=\nabla_{\partial_t}X=(X\ln\sigma)\partial_t$ and $\nabla_{\partial_t}\partial_t=-\sigma\grad_h\sigma+(\partial_t\ln\sigma)\partial_t$. Hence $\nabla_X\xi=0$ and $\nabla_\xi\xi=-\grad_h\ln\sigma=-\psi$, so that $\nabla_X\xi=-\eta(X)\psi$ for all $X\in TM$. Using $\nabla^hP=0$ we obtain, for $X,Y\in TN$,
        \[
            (\nabla_X\varphi)Y=0,\qquad (\nabla_X\varphi)\xi=0,\qquad (\nabla_\xi\varphi)Y=\omega(PY)\xi,\qquad (\nabla_\xi\varphi)\xi=\varphi\psi ,
        \]
        and this is exactly \eqref{eq:C}. Moreover $d\eta=d_N\sigma\wedge dt=\omega\wedge\eta$.

        \emph{Closedness.} $\omega=d_N\ln\sigma$ has no $dt$-component, so $d\omega=dt\wedge d_N(\partial_t\ln\sigma)$; this vanishes if and only if $\ln\sigma=A(t)+B(x)$. Then $\tau=\int e^{A}dt$ gives $g=h+e^{2B}d\tau^2$.

        \emph{Normality.} By the remark following Lemma~\ref{lem:basic}, the structure is normal if and only if $d\eta=0$, i.e.\ $\omega\wedge\eta=0$, i.e.\ $\omega=0$ (recall $\omega(\xi)=0$).
    \end{proof}

    \begin{remark}[The setting (Tw)]\label{rem:settingE}
        The splitting step of the proof only uses that $\xi$ is a unit vector field with $\nabla_X\xi=-\eta(X)\psi$, where $\eta=g(\xi,\cdot)$ and $\psi=-\nabla_\xi\xi$. Conversely, on every twisted product $g=h+\sigma^2dt^2$ over an \emph{arbitrary} pseudo-Riemannian manifold $(N,h)$, the field $\xi=\sigma^{-1}\partial_t$ has this property, with $\psi=\grad_h\ln\sigma$ and $\omega=d_N\ln\sigma$. We refer to a pseudo-Riemannian manifold with such a unit vector field as \emph{setting (Tw)}. Lemma~\ref{lem:basic}(i), (ii), (iv), Lemma~\ref{lem:Ricmodel}, Theorem~\ref{thm:compare}(i)--(iii), (v), Theorem~\ref{thm:nonclosed}, Theorem~\ref{thm:compact}(a), (b) and the results of Section~\ref{sec:solitons} are proved in setting (Tw). The almost (para)contact structure enters only through the (para-)K\"ahler structure of the base, i.e.\ through Theorem~\ref{thm:QE}, in Theorem~\ref{thm:nonclosedR}, Theorem~\ref{thm:compact}(c) and Corollaries~\ref{cor:R}, \ref{cor:P} and~\ref{cor:classR}.
    \end{remark}

    \begin{lemma}[Ricci curvature in the twisted model]\label{lem:Ricmodel}
        In the model \eqref{eq:model}, and more generally in setting (Tw), for $X,Y\in D$,
        \begin{gather*}
            S(X,Y)=\Ric_h(X,Y)-\frac{\Hess_h\sigma(X,Y)}{\sigma},\qquad S(\xi,D)=0,\\
            S(\xi,\xi)=-\dv\psi,\qquad \dv\psi=\frac{\Delta_h\sigma}{\sigma},
        \end{gather*}
        where $\Ric_h$, $\Hess_h$ and $\Delta_h$ are computed on the leaf $N\times\{t\}$.
    \end{lemma}

    \begin{proof}
        For $\sigma=\sigma(x)$ this is the classical warped-product formula, cf.\ \cite[Cor.~7.43]{ONe83}. In general, the leaves of $D$ are totally geodesic, so for $X,Y,Z\in D$ the curvature $R(X,Y)Z$ is that of $h$. Using $\nabla_X\xi=0$ for $X\in D$ one finds $R(Y,\xi)\xi=-\nabla_Y\psi-\omega(Y)\psi$, and hence $S(Y,Z)=\Ric_h(Y,Z)+g(R(Y,\xi)\xi,Z)=\Ric_h(Y,Z)-g(\nabla_Y\psi,Z)-\omega(Y)\omega(Z)$ for $Y,Z\in D$. With $\psi=\grad_h\ln\sigma$ we have $g(\nabla_Y\psi,Z)+\omega(Y)\omega(Z)=\sigma^{-1}\Hess_h\sigma(Y,Z)$. The remaining identities follow from Lemma~\ref{lem:basic}(iv) and $\dv\psi=\Delta_h\ln\sigma+g(\nabla_\xi\psi,\xi)=\Delta_h\ln\sigma+|\psi|^2=\Delta_h\sigma/\sigma$.
    \end{proof}

    \begin{remark}\label{rem:closed}
        Throughout, ``$\omega$ is closed'' refers to the specific $1$-form $\omega=g(\psi,\cdot)=-\nabla_\xi\eta$. In \cite[Def.~3.1]{Bel24} and \cite[Def.~2.1]{DB24} the definition of (almost) $C_{12}$-manifolds involves a closed $1$-form $\omega'$ with $d\eta=\omega'\wedge\eta$. If $\omega'$ is required to be $-\nabla_\xi\eta$, this is the restriction $\sigma=\alpha(t)a(x)$ of Theorem~\ref{thm:twisted}. If only the existence of some closed $\omega'$ is required, there is no restriction, since in the model \eqref{eq:model} the form $\omega'=\omega+\sigma^{-1}(\partial_t\ln\sigma)\,\eta$ is closed and satisfies $d\eta=\omega'\wedge\eta$.
    \end{remark}

    \begin{corollary}[Dimension three]\label{cor:dim3}
        Three-dimensional $C_{12}$- (resp.\ para-$C_{12}$-) manifolds are, locally, exactly the twisted products $h+\sigma^2dt^2$ with $h$ an arbitrary Riemannian (resp.\ Lorentzian) surface metric.
    \end{corollary}

    \begin{proof}
        In dimension three, $N$ is a surface, Riemannian in (R) and Lorentzian in (P). On a Riemannian surface every orthogonal almost complex structure is parallel. On a Lorentzian surface the two null line fields are parallel, since $h(\nabla_Xe_\pm,e_\pm)=0$ forces $\nabla_Xe_\pm\parallel e_\pm$; as $P=\pm\Id$ on them, $P$ is parallel. Hence every $(N,P,h)$ is K\"ahler, resp.\ para-K\"ahler, and Theorem~\ref{thm:twisted} applies.
    \end{proof}

    \begin{remark}\label{rem:Bel}
        In dimension three, Beldjilali \cite[Thm.~3.3 and Prop.~3.4]{Bel24} shows, using Olszak's formulas, that an almost contact metric structure is $C_{12}$ if and only if $\nabla_X\xi=-\eta(X)\psi$, if and only if $\nabla_{\varphi X}\xi=0$. Lemma~\ref{lem:Ricmodel} is consistent with \cite[(3.13)]{Bel24}: for $X,Y\in D$ it gives $S(X,Y)=(\frac r2+\dv\psi)g(X,Y)-\omega(X)\omega(Y)-g(\nabla_X\psi,Y)$, since the Gauss curvature of $h$ equals $\frac r2+\dv\psi$, and $S(\xi,\xi)=-\dv\psi$ follows from \cite[(3.13)]{Bel24} on using $e^{2\rho}=|\psi|^2$. For three-dimensional Walker manifolds of the class $\mathbb{G}_{12}$ see \cite{NB25}.
    \end{remark}
    \section{The $\eta$-reversed metric}\label{sec:compare}
    \noindent
    We call $\tilde g:=g-2\eta\otimes\eta$ the \emph{$\eta$-reversed metric}. It is the member $t=-2$ of the canonical variation $g_t=g+t\,\eta\otimes\eta$ of $g$ along the unit vector field $\xi$ studied by Olea \cite{Ole}. It coincides with $g$ on $D$, and since $g(\xi,\xi)=1$ we have $\tilde g(\xi,\xi)=-1$ and $\tilde g(\xi,\cdot)=-\eta$. Thus $\tilde g$ is Lorentzian in setting (R) and has signature $(n,n+1)$ in setting (P). Let $\nt$, $\tilde R$, $\tilde S$ and $\tilde r$ be the corresponding objects.

    \begin{theorem}\label{thm:compare}
        For every structure satisfying (C) (no closedness of $\omega$ is needed):
        \begin{itemize}
            \item[(i)] $\nt_XY=\nabla_XY+2\eta(X)\eta(Y)\psi$;
            \item[(ii)] $\tilde S=S+2\dv\psi\,\eta\otimes\eta$; in particular $\tilde S=S$ on $D\times D$;
            \item[(iii)] $\tilde r=r$;
            \item[(iv)] $\nt_X\xi=\eta(X)\psi$ and $(\nt_X\varphi)Y=\eta(X)\bigl\{\omega(\varphi Y)\xi-\eta(Y)\varphi\psi\bigr\}$;
            \item[(v)] $\tilde R(X,Y)Z=R(X,Y)Z-2\eta(Z)R(X,Y)\xi$.
        \end{itemize}
    \end{theorem}

    \begin{proof}
        (i) Let $Q:=\nt-\nabla$, a symmetric $(1,2)$-tensor, and $B(X,Y,Z):=\tilde g(Q(X,Y),Z)$. From $\nt\tilde g=0$ we get $(\nabla_X\tilde g)(Y,Z)=B(X,Y,Z)+B(X,Z,Y)$, hence
        \[
            B(X,Y,Z)=\tfrac12\bigl[(\nabla_X\tilde g)(Y,Z)+(\nabla_Y\tilde g)(X,Z)-(\nabla_Z\tilde g)(X,Y)\bigr].
        \]
        By Lemma~\ref{lem:basic}(i), $(\nabla_X\tilde g)(Y,Z)=-2\bigl[(\nabla_X\eta)(Y)\eta(Z)+\eta(Y)(\nabla_X\eta)(Z)\bigr]=2\eta(X)\bigl[\omega(Y)\eta(Z)+\eta(Y)\omega(Z)\bigr]$. Summing, $B(X,Y,Z)=2\eta(X)\eta(Y)\omega(Z)$. Since $\eta(\psi)=0$, $\tilde g(\psi,Z)=\omega(Z)$, and therefore $Q(X,Y)=2\eta(X)\eta(Y)\psi$.

        (ii) We have
        \begin{multline*}
            \tilde R(X,Y)Z=R(X,Y)Z+(\nabla_XQ)(Y,Z)-(\nabla_YQ)(X,Z)\\
            +Q(X,Q(Y,Z))-Q(Y,Q(X,Z)).
        \end{multline*}
        Since $\eta\circ Q=0$ (because $\eta(\psi)=0$), the last two terms vanish. By Lemma~\ref{lem:basic}(i),(iv),
        \begin{align*}
            \tr_X(\nabla_XQ)(Y,Z)&=2\eta(Y)\eta(Z)\dv\psi+2\bigl[(\nabla_\psi\eta)(Y)\eta(Z)+\eta(Y)(\nabla_\psi\eta)(Z)\bigr]\\
            &=2\eta(Y)\eta(Z)\dv\psi,\\
            \tr_X(\nabla_YQ)(X,Z)&=2\bigl[(\nabla_Y\eta)(\psi)+\eta(\nabla_Y\psi)\bigr]\eta(Z)\\
            &=2\bigl[-\eta(Y)|\psi|^2+\eta(Y)|\psi|^2\bigr]\eta(Z)=0,
        \end{align*}
        because $\nabla_\psi\eta=-\eta(\psi)\omega=0$. Thus $\tilde S=S+2\dv\psi\,\eta\otimes\eta$.

        (iii) $\tilde g^{-1}=g^{-1}-2\xi\otimes\xi$, so $\tilde r=r-2S(\xi,\xi)+2\dv\psi\,\tilde g^{-1}(\eta,\eta)=r+2\dv\psi-2\dv\psi=r$.

        (iv) $\nt_X\xi=-\eta(X)\psi+2\eta(X)\psi$. Further, $(\nt_X\varphi)Y=(\nabla_X\varphi)Y+Q(X,\varphi Y)-\varphi Q(X,Y)=(\nabla_X\varphi)Y-2\eta(X)\eta(Y)\varphi\psi$, and we use Lemma~\ref{lem:basic}(iii).

        (v) In the formula for $\tilde R$ in (ii) the quadratic terms vanish, and by Lemma~\ref{lem:basic}(i)
        \begin{multline*}
            (\nabla_XQ)(Y,Z)-(\nabla_YQ)(X,Z)\\
            =2\eta(Z)\Bigl\{\bigl[\eta(Y)\omega(X)-\eta(X)\omega(Y)\bigr]\psi+\eta(Y)\nabla_X\psi-\eta(X)\nabla_Y\psi\Bigr\}
        \end{multline*}
        (the terms containing $\omega(Z)$ cancel). By the proof of Lemma~\ref{lem:basic}(iv), $R(X,Y)\xi=\bigl[\eta(X)\omega(Y)-\eta(Y)\omega(X)\bigr]\psi-\eta(Y)\nabla_X\psi+\eta(X)\nabla_Y\psi$, so the term in braces equals $-R(X,Y)\xi$. Taking the trace of (v) in $X$ recovers (ii).
    \end{proof}

    \begin{remark}\label{rem:general}
        (1) The proofs of (i), (ii), (iii) and (v) only use $\nabla_X\xi=-\eta(X)\psi$ for a unit vector field $\xi$, i.e.\ setting (Tw) of Remark~\ref{rem:settingE}. Hence these formulas hold on every twisted product $N\times_\sigma I$ with an arbitrary base metric $h$ (with $\xi=\sigma^{-1}\partial_t$ and $\tilde g=h-\sigma^2dt^2$), without any almost contact structure.
        (2) In the model \eqref{eq:model}, (i) can also be read off from the Christoffel symbols: $g=h+\sigma^2dt^2$ and $\tilde g=h-\sigma^2dt^2$ have the same covariant derivatives $\nabla_XY$ and $\nabla_X\partial_t=\nabla_{\partial_t}X=(X\ln\sigma)\partial_t$ for $t$-independent $X,Y\in TN$, while $\nabla_{\partial_t}\partial_t=-\sigma\grad_h\sigma+(\partial_t\ln\sigma)\partial_t$ and $\nt_{\partial_t}\partial_t=\sigma\grad_h\sigma+(\partial_t\ln\sigma)\partial_t$ differ by $2\sigma^2\psi=2\eta(\partial_t)^2\psi$. The coordinate-free proof above shows that only $\nabla\xi=-\eta\otimes\psi$ is used.
        (3) For an arbitrary unit vector field $\xi$, the difference tensor of the Levi-Civita connections of $g$ and of $g_t=g+t\,\eta\otimes\eta$ is expressed in \cite[Prop.~2.3]{Ole} in terms of $\Lie_\xi g$ and $d\eta$; for $t=-2$ and $\nabla\xi=-\eta\otimes\psi$ it reduces to (i). The closed form of the curvature identities (ii) and (v) relies on the latter condition. In the warped case $\sigma=\sigma(x)$, (ii) also follows from the Ricci formulas of warped products \cite[Cor.~7.43]{ONe83}, whose proof applies to one-dimensional fibres of either sign.
    \end{remark}

    \noindent
    The structure $(\varphi,\xi,\eta,\tilde g)$ is an almost contact Lorentzian structure in setting (R): $\varphi^2=-\Id+\eta\otimes\xi$, $\eta(\xi)=1$, $\tilde g(\varphi X,\varphi Y)=\tilde g(X,Y)+\eta(X)\eta(Y)$ and $\tilde g(\xi,\cdot)=-\eta$; in setting (P) it satisfies $\varphi^2=\Id-\eta\otimes\xi$, $\tilde g(\varphi X,\varphi Y)=-\tilde g(X,Y)-\eta(X)\eta(Y)$ and $\tilde g(\xi,\cdot)=-\eta$. Conversely, every such structure arises in this way from $g:=\tilde g+2\eta\otimes\eta$.

    \begin{definition}\label{def:lorentz}
        Let $(\varphi,\xi,\eta,\tilde g)$ be a structure as above and put $\tilde F(X,Y,Z):=\tilde g((\nt_X\varphi)Y,Z)$. The structure is called an \emph{$\eta$-reversed $C_{12}$-structure} in setting (R), and an \emph{$\eta$-reversed para-$C_{12}$-structure} in setting (P), if
        \begin{equation}\label{eq:Ctilde}
            \tilde F(X,Y,Z)=\eta(X)\bigl\{\eta(Y)\tilde F(\xi,\xi,Z)-\eta(Z)\tilde F(\xi,\xi,Y)\bigr\}.
        \end{equation}
    \end{definition}

    \begin{proposition}\label{prop:lorentz}
        Let $g$ and $\tilde g=g-2\eta\otimes\eta$ be as above. Then $(\varphi,\xi,\eta,g)$ satisfies (C) if and only if $(\varphi,\xi,\eta,\tilde g)$ satisfies \eqref{eq:Ctilde}. In this case $\nt_X\xi=-\eta(X)\tilde\psi$ with $\tilde\psi:=-\nt_\xi\xi=-\psi$, and $\tilde g(\tilde\psi,\cdot)=-\omega$.
    \end{proposition}

    \begin{proof}
        Assume (C). By Theorem~\ref{thm:compare}(iv), $\nt_\xi\xi=\psi$ and $(\nt_X\varphi)Y=\eta(X)\{\omega(\varphi Y)\xi-\eta(Y)\varphi\psi\}$. Since $\tilde g(\xi,\cdot)=-\eta$ and $\tilde g(\varphi\psi,Z)=g(\varphi\psi,Z)=-\omega(\varphi Z)$, we get $\tilde F(X,Y,Z)=\eta(X)\{\eta(Y)\omega(\varphi Z)-\eta(Z)\omega(\varphi Y)\}$ and $\tilde F(\xi,\xi,Z)=\omega(\varphi Z)$, which is \eqref{eq:Ctilde}.

        Conversely, assume \eqref{eq:Ctilde}. The proof of Lemma~\ref{lem:basic}(i) only uses that $\varphi$ is skew, that the metric is non-degenerate on $D=\varphi(TM)$ and that it has a constant non-zero value on $\xi$; applied to $\tilde g$ it gives $\nt_X\xi=-\eta(X)\tilde\psi$ with $\tilde\psi:=-\nt_\xi\xi\in D$, and, as in Lemma~\ref{lem:basic}(iii), $(\nt_X\varphi)Y=\eta(X)\{\eta(Y)\varphi\tilde\psi-\tilde\omega(\varphi Y)\xi\}$ with $\tilde\omega:=\tilde g(\tilde\psi,\cdot)$. As $\eta=-\tilde g(\xi,\cdot)$, we have $(\nt_X\eta)Y=\eta(X)\tilde\omega(Y)$, and the computation in the proof of Theorem~\ref{thm:compare}(i), applied to $g=\tilde g+2\eta\otimes\eta$, gives $\nabla_XY=\nt_XY+2\eta(X)\eta(Y)\tilde\psi$. Hence $\psi=-\tilde\psi$, $\omega=-\tilde\omega$ and $(\nabla_X\varphi)Y=(\nt_X\varphi)Y-2\eta(X)\eta(Y)\varphi\tilde\psi=\eta(X)\{\omega(\varphi Y)\xi+\eta(Y)\varphi\psi\}$, which gives \eqref{eq:C} because $F(\xi,\xi,Z)=-\omega(\varphi Z)$.
    \end{proof}

    \noindent
    In setting (R), $\tilde g$ is Lorentzian, and an $\eta$-reversed $C_{12}$-structure may also be called a Lorentzian $C_{12}$-structure; in setting (P), $\tilde g$ is Lorentzian only for $n=1$, which is why we use the neutral term. In \cite[Def.~3.2, Prop.~3.3]{DB24} a ``Lorentz $C_{12}$-manifold'' is defined by $d\eta=\omega\wedge\eta$, $d\Phi=0$ and $N_\varphi=0$; Proposition~\ref{prop:lorentz} is the corresponding statement for the Chinea--Gonz\'alez type condition \eqref{eq:Ctilde}. Note that $\tilde\psi=-\psi$: the acceleration of $\xi$ changes sign.
    \section{$\eta$-Einstein structures}\label{sec:etaE}
    \noindent
    Throughout this section we use the model \eqref{eq:model}, with $\psi=\grad_h\ln\sigma$; the results up to and including Theorem~\ref{thm:compact}(a), (b) hold in setting (Tw) of Remark~\ref{rem:settingE}, whereas Theorem~\ref{thm:compact}(c) uses the K\"ahler structure of the base.

    \begin{definition}\label{def:etaE}
        A $C_{12}$- or para-$C_{12}$-manifold is called \emph{$\eta$-Einstein with constant $\lambda$} if there is a constant $\lambda\in\R$ such that $S(X,Y)=\lambda g(X,Y)$ for all $X,Y\in D$. By Lemma~\ref{lem:basic}(iv) this is equivalent to
        \begin{equation}\label{eq:etaE}
            S=\lambda g-\Theta\,\eta\otimes\eta,\qquad \Theta:=\lambda+\dv\psi ,
        \end{equation}
        and by Lemma~\ref{lem:Ricmodel} it is equivalent to the base equation
        \begin{equation}\label{eq:base}
            \Ric_h-\frac1\sigma\Hess_h\sigma=\lambda h\qquad\text{for every }t .
        \end{equation}
        The coefficient $\Theta$ of $\eta\otimes\eta$ is in general a function on $M$.
    \end{definition}

    \noindent
    Equation \eqref{eq:base} is \eqref{eq:QE} with $m=1$ and $u=\sigma(\cdot,t)$. For $m=1$, \eqref{eq:QE} is linear in $u$: it reads $\Hess_hu=u\,\mathcal R$ with $\mathcal R:=\Ric_h-\lambda h$, and for each $t$ the function $\sigma_t:=\sigma(\cdot,t)$ belongs to the vector space $\mathcal W$ of its solutions.

    \begin{remark}[Constancy of $\lambda$]\label{rem:nonconst}
        The constancy of $\lambda$ is essential, since it is needed in Lemma~\ref{lem:KK}; it is not automatic. For $\eta$-Einstein Sasakian manifolds of dimension at least five the coefficients are constant by Okumura's theorem \cite{Oku62}, but $C_{12}$-manifolds satisfy $\eta\wedge d\eta=0$, and there is no Schur-type lemma: over the flat K\"ahler manifold $\R^{2n}$ let $\sigma=1+|x|^2$. Then $\Hess_h\sigma=2h$, so by Lemma~\ref{lem:Ricmodel}
        \[
            S=-\frac2\sigma\,g+\frac{2-4n}{\sigma}\,\eta\otimes\eta ,
        \]
        i.e.\ $S(X,Y)=\lambda g(X,Y)$ on $D$ with the non-constant function $\lambda=-2/\sigma$, in every dimension $2n+1\geq3$. The results of this section do not apply to such structures.
    \end{remark}

    \begin{proposition}\label{prop:constcoef}
        Let $M$ be a $C_{12}$- or para-$C_{12}$-manifold with $S=\alpha_1g+\alpha_2\,\eta\otimes\eta$ for constants $\alpha_1,\alpha_2$. Then $\alpha_2=0$ or $\omega=0$. Thus an $\eta$-Einstein $C_{12}$- or para-$C_{12}$-manifold with constant coefficients is either Einstein, or it is normal and locally $g=h+dt^2$ with $\Ric_h=\alpha_1h$ and $\alpha_2=-\alpha_1$.
    \end{proposition}

    \begin{proof}
        Fix $t$. By Lemma~\ref{lem:Ricmodel}, $\sigma_t$ satisfies \eqref{eq:QE} on $N$ with $m=1$ and $\lambda=\alpha_1$, and $\dv\psi=\Delta_h\sigma_t/\sigma_t$. Since $S(\xi,\xi)=-\dv\psi=\alpha_1+\alpha_2$, Lemma~\ref{lem:KK} shows that $\sigma_t^2(\alpha_1+\dv\psi)=-\alpha_2\sigma_t^2$ is constant on $N$. If $\alpha_2\neq0$, then $\sigma_t$ is constant on $N$ for every $t$, so $\psi=\grad_h\ln\sigma=0$, i.e.\ $\omega=0$, and the structure is normal by Theorem~\ref{thm:twisted}. Reparametrising $t$ we get $\sigma=1$; then $S=\Ric_h$ on $D$ and $S(\xi,\xi)=0$, whence $\Ric_h=\alpha_1h$ and $\alpha_2=-\alpha_1$.
    \end{proof}

    \subsection{A global dichotomy}
    \begin{theorem}\label{thm:nonclosed}
        Let $(M^{2n+1},\varphi,\xi,\eta,g)$ be a connected $C_{12}$- or para-$C_{12}$-manifold which is $\eta$-Einstein with constant $\lambda$, and let $\Theta=\lambda+\dv\psi$ be as in \eqref{eq:etaE}.
        \begin{itemize}
            \item[(i)] In every local model \eqref{eq:model}, $\Theta=\scal_h-(2n-1)\lambda$ is a function on $N$, and $\sigma(x,t)^2\,\Theta(x)$ depends only on $t$.
            \item[(ii)] $d\Theta=-2\Theta\,\omega$ on $M$.
            \item[(iii)] Either $\Theta\equiv0$, i.e.\ $S=\lambda g$ on $M$, or $\Theta$ has no zeros and $\omega=-\frac12\,d\ln|\Theta|$; in the latter case $\omega$ is exact on $M$.
            \item[(iv)] In particular, if $d\omega\neq0$ at some point of $M$, then $S=\lambda g$ on all of $M$.
        \end{itemize}
        The same holds in setting (Tw), i.e.\ for every twisted product $h+\sigma^2dt^2$ over an arbitrary pseudo-Riemannian manifold $(N,h)$ that satisfies \eqref{eq:base}.
    \end{theorem}

    \begin{proof}
        (i) Taking the $h$-trace of \eqref{eq:base} gives, by Lemma~\ref{lem:Ricmodel}, $\dv\psi=\Delta_h\sigma/\sigma=\scal_h-2n\lambda$, so $\Theta=\scal_h-(2n-1)\lambda$; the scalar curvature $\scal_h$ of $h$ does not depend on $t$. For fixed $t$, Lemma~\ref{lem:KK} with $m=1$ and $u=\sigma_t$ shows that $\mu=\sigma\Delta_h\sigma+\lambda\sigma^2=\sigma^2(\dv\psi+\lambda)=\sigma^2\Theta$ is constant on $N$, i.e.\ $\sigma(x,t)^2\Theta(x)=\mu(t)$.

        (ii) By (i), $\xi(\Theta)=0=-2\Theta\,\omega(\xi)$ and $d_N\Theta=d_N(\mu(t)\sigma^{-2})=-2\Theta\,d_N\ln\sigma=-2\Theta\,\omega$ in every local model. As $\Theta$ and $\omega$ are globally defined, (ii) holds on $M$.

        (iii) Along every smooth curve $\gamma$ in $M$, the function $y:=\Theta\circ\gamma$ satisfies the linear equation $y'=-2\,\omega(\gamma')\,y$; hence $y\equiv0$ if $y$ vanishes at one point. Therefore the zero set of $\Theta$ is open and closed, and either $\Theta\equiv0$ or $\Theta$ has no zeros. In the first case $S=\lambda g$ by \eqref{eq:etaE}; in the second, (ii) gives $d\ln|\Theta|=-2\omega$.

        (iv) If $\Theta$ had no zeros, $\omega$ would be closed. Only \eqref{eq:base} and Lemmas~\ref{lem:Ricmodel} and~\ref{lem:KK} were used, and the proof of Lemma~\ref{lem:Ricmodel} applies in setting (Tw); this gives the last assertion.
    \end{proof}

    \begin{remark}[Warped products and static potentials]\label{rem:HPW}
        (1) For every $v\in\mathcal W$ the function $\mu(v)=v\Delta_hv+\lambda v^2$ is constant by Lemma~\ref{lem:KK}, and, since $\Delta_hv=v(\scal_h-2n\lambda)$, one has the pointwise identity $\mu(v)=v^2\bigl(\scal_h-(2n-1)\lambda\bigr)$. The quadratic form $\mu$ on the space $\mathcal W$ of solutions of \eqref{eq:QE} was studied by He, Petersen and Wylie for Riemannian metrics \cite{HPW15b} (see also \cite{HPW14}); the degenerate case $m=1$, including its relation to static metrics, is discussed in \cite[Appendix~2]{HPW15b}. Theorem~\ref{thm:nonclosed}(i) is an instance of this structure: where $\scal_h\neq(2n-1)\lambda$, all elements of $\mathcal W$ are proportional to $|\scal_h-(2n-1)\lambda|^{-1/2}$, so that $\sigma$ separates variables; a non-closed $\omega$ requires $\dim\mathcal W\geq2$ and forces $\scal_h=(2n-1)\lambda$.
        (2) Let $L_h^*v:=-(\Delta_hv)h+\Hess_hv-v\Ric_h$ be the formal adjoint of the linearised scalar curvature operator \cite{FM75}. For $v\in\mathcal W$ one finds $L_h^*v=-v\bigl(\scal_h-(2n-1)\lambda\bigr)h$. Hence, where $\Theta=0$, every $\sigma_t$ is a static potential of $h$ in the sense of Fischer--Marsden and Corvino \cite{Cor00,FM75}, and conversely every static potential of a metric with constant scalar curvature $(2n-1)\lambda$ belongs to $\mathcal W$. In particular, near a point where $d\omega\neq0$, the base $(N,h)$ is a static space whose static potentials form a space of dimension at least two.
        (3) Compared with \cite{HPW14,HPW15b}, the new points are that the argument is carried out in arbitrary signature and yields the global identity (ii) on $M$, and that the (para-)K\"ahler structure of the base gives the complete local description in Theorem~\ref{thm:nonclosedR}(a), where the factor carrying the non-trivial solutions is two-dimensional; the para-K\"ahler case, in particular its isotropic part, has no counterpart there.
    \end{remark}

    \subsection{Non-closed $\omega$}
    \noindent
    Where $\omega$ is not closed, the base equation has several independent solutions $\sigma_t$, and Theorem~\ref{thm:QE} gives the local structure.

    \begin{theorem}\label{thm:nonclosedR}
        Let $M$ be as in Theorem~\ref{thm:nonclosed} and let $W\subset M$ be the open set where $d\omega\neq0$; if $W\neq\emptyset$, then $S=\lambda g$ on $M$ by Theorem~\ref{thm:nonclosed}(iv).
        \begin{itemize}
            \item[(a)] In setting (R), every point of a dense open subset of $W$ has a neighbourhood isomorphic, as a $C_{12}$-manifold, to the product of a $C_{12}$-manifold $M_1$ and a K\"ahler manifold $(K,h_K)$, possibly of dimension zero, where either
            \begin{itemize}
                \item[(a1)] $\lambda\neq0$, $M_1=(\Sigma\times I,\,h_\Sigma+\sigma^2dt^2)$ with $\Sigma$ a surface of constant curvature $\lambda/2$, $\sigma(\cdot,t)$ belongs to the three-dimensional space $V_\Sigma:=\{v:\ \Hess_{h_\Sigma}v=-\frac\lambda2\,v\,h_\Sigma\}$ for all $t$, and $\Ric_K=\lambda h_K$; or
                \item[(a2)] $\lambda=0$, $M_1=\bigl(\R^{2k}\times I,\ \sum_i(dx^i)^2+(a(t)+\sum_ib_i(t)x^i)^2dt^2\bigr)$ with $1\leq k\leq n$, and $\Ric_K=0$.
            \end{itemize}
            In both cases $M_1$ has constant curvature $\lambda/2$. In (a2), $2k$ is the rank of the distribution spanned by the parallel vector fields $\nabla\sigma_t$, $P\nabla\sigma_t$ ($t\in I$), which is constant on every connected open set of the base where $\Ric_h-\lambda h=0$; every value $1\leq k\leq n$ occurs (Example~\ref{ex:flatnc}).
            \item[(b)] In setting (P), let $p=(x_0,t_0)\in W$. If $\mathcal R\neq0$ at $x_0$ and $h(\nabla\sigma_t,\nabla\sigma_t)(x_0)\neq0$ for some $t$, then near $p$ the conclusion of (a1) holds, with $\Sigma$ a Lorentzian surface and $K$ para-K\"ahler. If $\mathcal R=0$ near $x_0$, or if $\nabla\sigma_t$ is null near $x_0$ for all $t$, then $\lambda=0$ and $g$ is Ricci-flat. The points of $W$ of these three kinds form a dense open subset of $W$. In the last two cases we do not obtain a normal form; in this sense (b) is a partial classification (Example~\ref{ex:ncP}).
        \end{itemize}
    \end{theorem}

    \begin{proof}
        Let $p=(x_0,t_0)\in W$ and let $U\ni x_0$ be a connected neighbourhood. The equation $\Hess_hv=v\mathcal R$ is a linear system whose solutions on $U$ are determined by $v(x_0)$ and $dv(x_0)$, and by \eqref{eq:base} it is satisfied by $\sigma_t$ for every $t$. If $d\sigma_t(x_0)=0$ for all $t$ near $t_0$, then also $\partial_td\sigma_t(x_0)=0$, so $d_N(\partial_t\ln\sigma)$ vanishes at $p$, contradicting $d\omega=dt\wedge d_N(\partial_t\ln\sigma)\neq0$. Hence $d\sigma_{t_1}(x_0)\neq0$ for some $t_1$, which we fix. By Lemma~\ref{lem:Ricform}, $\mathcal R(P\cdot,P\cdot)=-\eps\mathcal R$, so that $\mathcal R$, viewed as an endomorphism, commutes with $P$.

        (a) Suppose first that $\mathcal R\neq0$ at $x_0$. As $h$ is Riemannian, $\nabla\sigma_{t_1}$ is not null near $x_0$, and Theorem~\ref{thm:QE}(ii), applied with $m=1$ to $u=\sigma_{t_1}$, gives a parallel plane field $\mathcal Z$ with $A(\mathcal Z^\perp)=0$, where $A=\nabla\nabla\sigma_{t_1}=\sigma_{t_1}\mathcal R$; hence $\operatorname{Im}\mathcal R\subset\mathcal Z$. As $\operatorname{Im}\mathcal R$ is $P$-invariant and $P$ has no real eigenvectors, $\operatorname{Im}\mathcal R=\mathcal Z$ near $x_0$. By the local de Rham decomposition, which holds in every signature \cite{PR93}, $N=\Sigma\times K$ locally, with $T\Sigma=\mathcal Z=\operatorname{Im}\mathcal R$ and $TK=\ker\mathcal R$. Let $v$ solve $\Hess_hv=v\mathcal R$. For $X\in TK$ we have $\nabla_X\nabla v=v\,\mathcal R(X)=0$, and for $X\in T\Sigma$ we have $\nabla_X\nabla v\in T\Sigma$; hence the $TK$-component of $\nabla v$ is parallel, $\Hess_hv$ has no mixed terms, and $v=v_1(\Sigma)+v_2(K)$ with $\Hess_Kv_2=0$. Restricting $\Hess_hv=v\mathcal R$ to $T\Sigma$ gives $\Hess_\Sigma v_1=(v_1+v_2)\,\mathcal R_\Sigma$ with $\mathcal R_\Sigma:=\mathcal R|_{T\Sigma}\neq0$, so $v_2$ is constant. Thus $\sigma=\sigma(x_\Sigma,t)$, $g=(h_\Sigma+\sigma^2dt^2)+h_K$, $\Ric_K=\lambda h_K$ because $\mathcal R$ vanishes on $TK$, and $M_1:=(\Sigma\times I,h_\Sigma+\sigma^2dt^2)$ is an $\eta$-Einstein $C_{12}$-manifold with the same $\lambda$ and $d\omega\neq0$ near $p$. Theorem~\ref{thm:nonclosed}(i),(iv), applied to $M_1$, gives $\scal_{h_\Sigma}=\lambda$ near the projection of $p$, i.e.\ $\Sigma$ has constant curvature $\lambda/2$ there; hence $\mathcal R_\Sigma=-\frac\lambda2h_\Sigma$, $\lambda\neq0$, $\sigma(\cdot,t)\in V_\Sigma$, and $M_1$ is Einstein with $S=\lambda g$, hence of constant curvature $\lambda/2$. On a surface of constant curvature $\lambda/2\neq0$ the space $V_\Sigma$ is three-dimensional. This is (a1).

        Suppose now that $\mathcal R=0$ near $x_0$. Then $\Hess_h\sigma_t=0$ for all $t$, and $\nabla\sigma_{t_1}$ is a non-zero parallel vector field, so $\Ric_h(\nabla\sigma_{t_1},\cdot)=0$; as $\Ric_h=\lambda h$, we get $\lambda=0$. The parallel vector fields $\nabla\sigma_t$ and $P\nabla\sigma_t$ ($t\in I$) span a parallel, flat, $P$-invariant distribution $E_0$; since parallel vector fields are determined by their values at one point, linear relations among finitely many of them hold everywhere if they hold at one point, so the rank of $E_0$ is constant on connected sets, say $2k$ with $1\leq k\leq n$, and $E_0$ is non-degenerate in setting (R). By the local de Rham decomposition, $N=\R^{2k}\times K$ locally, with $T\R^{2k}=E_0$ and $\Ric_K=0$, and each $\sigma_t$ is an affine function on $\R^{2k}$. Finally, $M_1=(\R^{2k}\times I,\sum_i(dx^i)^2+\sigma^2dt^2)$ is flat: for $X,Y,Z\in D$, $R(X,Y)Z$ is the curvature of the flat leaves, $R(X,Y)\xi=0$ by the proof of Lemma~\ref{lem:basic}(iv), and $R(X,\xi)\xi=-\sigma^{-1}(\Hess_h\sigma(X,\cdot))^\sharp=0$ by the proof of Lemma~\ref{lem:Ricmodel}; the remaining components vanish by the symmetries of $R$. This is (a2). Since $\{\mathcal R\neq0\}$ and the interior of $\{\mathcal R=0\}$ form a dense open subset of $N$, (a) follows.

        (b) Suppose that $\mathcal R\neq0$ at $x_0$ and that $\nabla\sigma_{t_1}$ is not null at $x_0$ (choosing $t_1$ accordingly). As in (a), Theorem~\ref{thm:QE}(ii) gives $\operatorname{Im}\mathcal R\subset\mathcal Z$, where now $\mathcal Z$ is a non-degenerate $P$-invariant plane. If $\operatorname{Im}\mathcal R$ were a line, it would be $P$-invariant and $\mathcal R=\theta\,\alpha\otimes\alpha$ for a function $\theta\neq0$ and a $1$-form $\alpha$; then $\mathcal R(P\cdot,P\cdot)=-\mathcal R$ would give $\alpha(PX)^2=-\alpha(X)^2$ for all $X$, i.e.\ $\alpha=0$. Hence $\operatorname{Im}\mathcal R=\mathcal Z$, and the argument for (a1) applies verbatim, $\Sigma$ being a Lorentzian surface. If $\mathcal R=0$ near $x_0$, then $\lambda=0$ as in (a). If $\nabla\sigma_t$ is null near $x_0$ for all $t$, then Theorem~\ref{thm:QE}(iii), applied to $\sigma_{t_1}$ with $m=1$, gives $\lambda=0$. In both cases $S=\lambda g=0$ by Theorem~\ref{thm:nonclosed}(iv). Finally, let $A:=\{\mathcal R\neq0\}$ and let $G$ be the open set of points $(x,t')$ of $W$ such that $h(\nabla\sigma_t,\nabla\sigma_t)(x)\neq0$ for some $t$. If an open set $O\subset W$ does not meet $A\cap G$, then either $O\cap G\neq\emptyset$, and $O\cap G$ is an open subset of $\{\mathcal R=0\}$, or $O\cap G=\emptyset$, and $O$ lies in the interior of $W\setminus G$; this proves the density statement.
    \end{proof}

    \subsection{Compact manifolds}
    \begin{theorem}\label{thm:compact}
        Let $(M,\varphi,\xi,\eta,g)$ be a compact connected $C_{12}$- or para-$C_{12}$-manifold which is $\eta$-Einstein with constant $\lambda$, and let $\Theta=\lambda+\dv\psi$.
        \begin{itemize}
            \item[(a)] If $\lambda=0$, then $g$ is Ricci-flat and $\dv\psi=0$.
            \item[(b)] If $\lambda\neq0$, then $\Theta$ has no zeros and has the sign of $\lambda$; in particular $g$ is not Einstein. Moreover $\omega=df$ with $f:=-\frac12\ln|\Theta|$, and $e^f\xi$ is a nowhere vanishing Killing vector field of $g$ and of $\tilde g$.
            \item[(c)] In setting (R), $\omega=0$ for every $\lambda$. Hence the structure is normal, $\nabla\xi=0$ and $\nabla\varphi=0$, $S=\lambda(g-\eta\otimes\eta)$, and locally $g=h+dt^2$ with $\Ric_h=\lambda h$.
        \end{itemize}
        In particular, a compact Einstein $C_{12}$- or para-$C_{12}$-manifold is Ricci-flat, and a compact $\eta$-Einstein $C_{12}$-manifold with constant $\lambda$ is Einstein only if it is Ricci-flat.
    \end{theorem}

    \begin{proof}
        By the divergence theorem, $\int_M\Theta=\int_M(\lambda+\dv\psi)=\lambda\operatorname{vol}(M)$. If $\lambda=0$ and $\Theta$ had no zeros, $\Theta$ would have a constant sign and $\int_M\Theta\neq0$; hence $\Theta\equiv0$ by Theorem~\ref{thm:nonclosed}(iii), and $S=0$, $\dv\psi=\Theta-\lambda=0$. If $\lambda\neq0$, then $\Theta\equiv0$ is impossible, so $\Theta$ has no zeros, its sign is that of $\int_M\Theta=\lambda\operatorname{vol}(M)$, and $S(\xi,\xi)=\lambda-\Theta\neq\lambda$. By Theorem~\ref{thm:nonclosed}(iii), $\omega=df$. Since $\Lie_\xi g=-\eta\otimes\omega-\omega\otimes\eta$ by Lemma~\ref{lem:basic}(i), we get $\Lie_{e^f\xi}g=e^f\bigl(-\eta\otimes\omega-\omega\otimes\eta+df\otimes\eta+\eta\otimes df\bigr)=0$; moreover $\Lie_{e^f\xi}\eta=d(e^f)+e^f\iota_\xi(\omega\wedge\eta)=e^f(df-\omega)=0$, so $e^f\xi$ is also a Killing field of $\tilde g=g-2\eta\otimes\eta$. This proves (a) and (b).

        (c) Now $g$ is Riemannian. We first collect three facts about the leaves of $D=\ker\eta$, which carry the induced Riemannian metric $h$ and the K\"ahler structure $P=\varphi|_D$ (Theorem~\ref{thm:twisted}).

        \emph{(1) The leaves are complete.} A leaf $L$ is totally geodesic (proof of Theorem~\ref{thm:twisted}), so a geodesic of $M$ which starts tangent to $L$ stays in $L$ and is a geodesic of $L$: in every foliation chart the plaques are totally geodesic, and the uniqueness of geodesics shows that the geodesic cannot leave the plaque through its starting point. Since $M$ is compact, geodesics of $M$ are defined for all time; hence $(L,h)$ is geodesically complete, and so is its universal cover.

        \emph{(2) Solutions of the base equation on the leaves.} On a leaf $L$ the form $\omega|_{TL}$ is closed, since $d\omega=dt\wedge d_N(\partial_t\ln\sigma)$ vanishes on $TL$ in every local model, and on a plaque $N\times\{t\}$ it equals $d\ln\sigma_t$ with $\sigma_t=\sigma(\cdot,t)$ (Theorem~\ref{thm:twisted}). If $v>0$ is a function on a connected open subset of $L$, or of its universal cover, with $d\ln v=\omega|_{TL}$, then $v$ is a positive constant multiple of $\sigma_t$ on every plaque; as \eqref{eq:base} is linear in $\sigma_t$, $v$ satisfies $\Hess_hv=v\,(\Ric_h-\lambda h)$, and $\Delta_hv=v\,\dv\psi$ by Lemma~\ref{lem:Ricmodel}.

        \emph{(3) The Hessian of such solutions.} Let $v>0$ satisfy $\Hess_hv=v\,(\Ric_h-\lambda h)$ on a connected K\"ahler manifold $(L,h,P)$ with $h$ Riemannian, and put $\zeta:=\nabla v$, $A:=\nabla\zeta$. Wherever $\zeta\neq0$ we have $h(\zeta,\zeta)>0$, and Theorem~\ref{thm:QE}(i), (ii) with $m=1$ gives $AP=PA$, $A(\mathcal Z)\subset\mathcal Z$ and $A(\mathcal Z^\perp)=0$ for $\mathcal Z=\Span\{\zeta,P\zeta\}$. On the plane $\mathcal Z$, $P$ is a rotation by a right angle, and $A|_{\mathcal Z}$ is symmetric and commutes with $P$; hence $A|_{\mathcal Z}=\kappa\,\Id$, i.e.
        \begin{equation}\label{eq:Akappa}
            A=\kappa\,\pi_{\mathcal Z}\quad\text{on }\{\zeta\neq0\},\qquad \kappa:=\tfrac12\Delta_hv ,
        \end{equation}
        where $\pi_{\mathcal Z}$ is the orthogonal projection onto $\mathcal Z$. Moreover $V:=P\zeta$ is a Killing field of $(L,h)$ with $\nabla V=PA$, by Theorem~\ref{thm:QE}(i).

        \emph{Case $\lambda<0$.} By (b), $\Theta=-e^{-2f}$. As $\psi=\nabla f$, we have $\Delta_gf=\dv\psi=\Theta-\lambda=|\lambda|-e^{-2f}$. At a maximum point of $f$ we have $\Delta_gf\leq0$, hence $e^{-2\max f}\geq|\lambda|$; at a minimum point $\Delta_gf\geq0$, hence $e^{-2\min f}\leq|\lambda|$. Thus $\min f\geq\max f$, $f$ is constant and $\omega=df=0$.

        \emph{Case $\lambda>0$.} By (b), $\Theta>0$ and $\omega=df$ with $f=-\frac12\ln\Theta$. Fix a leaf $L$ and put $u:=e^f|_L=\Theta^{-1/2}|_L$, a smooth function on $L$ with $d\ln u=\omega|_{TL}$ and $0<\min_Me^f\leq u\leq\max_Me^f$. By Theorem~\ref{thm:nonclosed}(i), $\sigma_t^2\Theta=\mu(t)>0$ on the plaque with parameter $t$, so $u=\sigma_t/\sqrt{\mu(t)}$ there, and by (2)
        \begin{equation}\label{eq:leafQE}
            \Ric_h-\tfrac1u\Hess_hu=\lambda h,\qquad u\,\Delta_hu+\lambda u^2=u^2\,(\dv\psi+\lambda)=u^2\Theta=1\qquad\text{on }L .
        \end{equation}

        \emph{Step 1: $L$ is compact.} Consider the warped product $E:=L\times\R$ with the Riemannian metric $g_E:=h+u^2dt^2$. By Lemma~\ref{lem:Ricmodel} and \eqref{eq:leafQE}, its Ricci tensor satisfies $S_E(X,X)=\Ric_h(X,X)-\frac1u\Hess_hu(X,X)=\lambda\,h(X,X)$ for $X\in TL$. Let $p,q\in L$ and let $\gamma\colon[0,\ell]\to L$ be a unit-speed minimizing geodesic from $p$ to $q$, which exists by (1) and the Hopf--Rinow theorem. The curve $\tilde\gamma:=(\gamma,0)$ is a geodesic of $E$, because $L\times\{0\}$ is totally geodesic in the warped product \cite[Ch.~7]{ONe83}, and it is minimizing in $E$, because the projection $E\to L$ does not increase lengths. Suppose that $\ell>\pi\sqrt{2n/\lambda}$. For parallel orthonormal fields $E_1,\dots,E_{2n}$ along $\tilde\gamma$ orthogonal to $\tilde\gamma'$ and the proper variation fields $V_i:=\sin(\pi s/\ell)E_i$, the index forms add up to
        \[
            \sum_{i=1}^{2n}I(V_i,V_i)=\int_0^\ell\Bigl(2n\,\frac{\pi^2}{\ell^2}\cos^2\frac{\pi s}{\ell}-\sin^2\frac{\pi s}{\ell}\,S_E(\tilde\gamma',\tilde\gamma')\Bigr)ds=\frac{2n\pi^2-\lambda\ell^2}{2\ell}<0 ,
        \]
        so some $V_i$ yields curves from $(p,0)$ to $(q,0)$ shorter than $\tilde\gamma$, a contradiction. Hence $\operatorname{diam}(L,h)\leq\pi\sqrt{2n/\lambda}$, and $L$ is compact. (This is the case $m=1$ of a Myers-type theorem of Qian \cite{Qian97} for the tensor $\Ric_h-\frac mu\Hess_hu$; the warped product reduces it to the classical argument.)

        \emph{Step 2: a first integral.} Put $\zeta:=\nabla u$ and $A:=\nabla\zeta$ on $L$. On $\{\zeta\neq0\}$ we have $A\zeta=\kappa\zeta$ by \eqref{eq:Akappa}, hence $d\,h(\zeta,\zeta)=2h(A\zeta,\cdot)=2\kappa\,du=\Delta_hu\,du=(u^{-1}-\lambda u)\,du$ by \eqref{eq:leafQE}. Thus the differential of the smooth function $h(\zeta,\zeta)-\ln u+\frac\lambda2u^2$ vanishes on $\{\zeta\neq0\}$, and trivially on the interior of $\{\zeta=0\}$; as the union of these two sets is dense and $L$ is connected, there is a constant $C$ with
        \begin{equation}\label{eq:firstint}
            h(\nabla u,\nabla u)=\Phi(u),\qquad \Phi(x):=\ln x-\tfrac\lambda2x^2+C,\qquad\text{on }L .
        \end{equation}
        Note that $\Phi''(x)=-x^{-2}-\lambda<0$ and that, by \eqref{eq:leafQE}, $\kappa=\frac12\Delta_hu=\frac12\Phi'(u)$ on $L$.

        \emph{Step 3: the period of the Killing field $V=P\nabla u$.} Suppose that $u$ is not constant on $L$, and let $u_1:=\min_Lu<u_2:=\max_Lu$. Since $\nabla u=0$ at extremal points, \eqref{eq:firstint} gives $\Phi(u_1)=\Phi(u_2)=0$, and the strict concavity of $\Phi$ gives $\Phi'(u_1)>0>\Phi'(u_2)$. The closed set $\{u=u_1\}$ is a non-empty proper subset of the connected manifold $L$, hence not open, so it has a boundary point $p_1$; every neighbourhood of $p_1$ contains points with $\nabla u\neq0$, since otherwise $u$ would be constant near $p_1$. The identity $A\circ A=\kappa A$ holds on $\{\zeta\neq0\}$ by \eqref{eq:Akappa}, hence at $p_1$ by continuity, and $\tr A_{p_1}=\Delta_hu(p_1)=2\kappa(p_1)$ with $\kappa(p_1)=\frac12\Phi'(u_1)\neq0$. Therefore the symmetric endomorphism $A_{p_1}$ of $T_{p_1}L$ has the eigenvalue $\kappa(p_1)$ with multiplicity two and the eigenvalue $0$ with multiplicity $2n-2$, and its $\kappa(p_1)$-eigenspace $\mathcal Z_1$ is $P$-invariant because $AP=PA$. The Killing field $V$ vanishes at $p_1$, and $(\nabla V)_{p_1}=PA_{p_1}$ acts on $\mathcal Z_1$ as $\kappa(p_1)P$ and vanishes on $\mathcal Z_1^\perp$. Its flow $\phi_t$ is defined for all $t\in\R$ because $L$ is compact, it fixes $p_1$, and $d(\phi_t)_{p_1}=\exp\bigl(t(\nabla V)_{p_1}\bigr)$ (the linearisation of a flow at a zero of the vector field) is the rotation by the angle $\kappa(p_1)\,t$ in the plane $\mathcal Z_1$ and the identity on $\mathcal Z_1^\perp$. An isometry of the connected manifold $L$ which fixes a point and has differential the identity there is the identity; hence $\phi_t=\mathrm{id}_L$ if and only if $\kappa(p_1)\,t\in2\pi\mathbb Z$. The same argument at a boundary point $p_2$ of $\{u=u_2\}$ shows that $\phi_t=\mathrm{id}_L$ if and only if $\kappa(p_2)\,t\in2\pi\mathbb Z$, where $\kappa(p_2)=\frac12\Phi'(u_2)\neq0$. Comparing the two conditions gives $|\kappa(p_1)|=|\kappa(p_2)|$, i.e.\ $\Phi'(u_1)=-\Phi'(u_2)$, which reads $u_1^{-1}+u_2^{-1}=\lambda(u_1+u_2)$, i.e.\ $u_1u_2=1/\lambda$.

        \emph{Step 4: contradiction.} Write $u_1=\lambda^{-1/2}e^{-s}$ and $u_2=\lambda^{-1/2}e^{s}$ with $s>0$. Then $\Phi(u_1)=\Phi(u_2)$ reads $-s-\frac12e^{-2s}=s-\frac12e^{2s}$, i.e.\ $\sinh2s=2s$, which is impossible for $s>0$. Hence $u$ is constant on every leaf, so $\omega|_{TL}=d\ln u=0$ for every leaf $L$; as also $\omega(\xi)=0$, we get $\omega=0$.

        \emph{Case $\lambda=0$.} By (a), $\Theta\equiv0$ and $\dv\psi=0$. Let $L$ be a leaf and $\tilde L$ its universal cover with the pulled-back metric and K\"ahler structure, which is complete by (1). The pull-back of the closed form $\omega|_{TL}$ is exact on $\tilde L$, say equal to $d\ln v$ with $v>0$. By (2), $v$ satisfies $\Hess_hv=v\Ric_h$ and $\Delta_hv=v\dv\psi=0$, so $\kappa=0$ and $\Hess_hv=0$ on $\{dv\neq0\}$ by \eqref{eq:Akappa}; trivially $\Hess_hv=0$ on the interior of $\{dv=0\}$, hence $\Hess_hv=0$ on $\tilde L$ by density. For every geodesic $\gamma$ of $\tilde L$, $(v\circ\gamma)''=\Hess_hv(\gamma',\gamma')=0$, so $v\circ\gamma$ is an affine function on $\R$, which is constant because it is positive (cf.\ \cite[Remark~3.3(3)]{SS-QE}). Hence $v$ is constant, $\omega|_{TL}=0$ for every leaf, and $\omega=0$.

        The remaining assertions of (c) follow from Lemma~\ref{lem:basic}(i),(iii), Theorem~\ref{thm:twisted} and Lemma~\ref{lem:Ricmodel}: $\omega=0$ gives $\psi=0$, $\sigma=1$ after reparametrising $t$, $S=\Ric_h=\lambda h$ on $D$ and $S(\xi,\cdot)=0$. The last statements of the theorem follow from (a), (b) and (c).
    \end{proof}

    \noindent
    Compact examples exist for every $\lambda$, e.g.\ $K\times S^1$ with a compact K\"ahler--Einstein manifold $K$ with $\Ric_K=\lambda h_K$, and by (c) every compact $\eta$-Einstein $C_{12}$-manifold with constant $\lambda$ is locally of this form. Compactness cannot be weakened to finite volume: in Example~\ref{ex:S3dS}, $\lambda=2$ and $\omega\neq0$, and the structure lives on the complement of a circle in $S^3$ but does not extend to $S^3$, since $\xi=\partial_t/\cos r$ becomes singular where $\cos r=0$. The proof of (c) uses the Riemannian signature of $g$ in the maximum principle, in Myers' argument and in \eqref{eq:Akappa}; in setting (P) we do not know whether compact examples with $\omega\neq0$ exist (Section~\ref{sec:open}).

    \subsection{Closed $\omega$}
    \begin{corollary}[Setting (R)]\label{cor:R}
        Let $(M^{2n+1},\varphi,\xi,\eta,g)$ be a $C_{12}$-manifold with $d\omega=0$ and $n\geq2$ which is $\eta$-Einstein with constant $\lambda$. On every connected open set where $\omega\neq0$, $M$ is locally isomorphic, as a $C_{12}$-manifold, to the product of a three-dimensional $\eta$-Einstein $C_{12}$-manifold $(\Sigma\times I,\,h_\Sigma+a^2dt^2)$, $a=a(\Sigma)$, and a K\"ahler--Einstein manifold $(K,h_K)$ with $\Ric_K=\lambda h_K$. The same splitting holds for $\tilde g$. The three-dimensional factor is Einstein if and only if $\mu=a\Delta_ha+\lambda a^2$ vanishes.
    \end{corollary}

    \begin{proof}
        By Theorem~\ref{thm:twisted}, $d\omega=0$ gives $\sigma=a(x)$ after reparametrising $t$, and \eqref{eq:base} is \eqref{eq:QE} with $m=1$ and $u=a$ on a K\"ahler base. As $h$ is positive definite, $h(\nabla a,\nabla a)>0$ wherever $da\neq0$, i.e.\ wherever $\omega=d_N\ln a\neq0$. Theorem~\ref{thm:QE}(ii) gives $N=\Sigma\times K$ locally with $a=a(\Sigma)$; hence $g=(h_\Sigma+a^2dt^2)+h_K$ and $\tilde g=(h_\Sigma-a^2dt^2)+h_K$, and $\varphi$ is the product of the $C_{12}$-structure on $\Sigma\times I$ given by the converse part of Theorem~\ref{thm:twisted} and the K\"ahler structure $P_K$ of $K$. The factor $\Sigma\times I$ is $\eta$-Einstein with the same $\lambda$, and it is Einstein if and only if $\lambda+\dv\psi=\lambda+\Delta_ha/a=0$, i.e.\ if and only if $\mu=0$.
    \end{proof}

    \begin{corollary}[Setting (P)]\label{cor:P}
        Let $(M^{2n+1},\varphi,\xi,\eta,g)$ be a para-$C_{12}$-manifold with $d\omega=0$ which is $\eta$-Einstein with constant $\lambda$, and let $W$ be a connected open set on which $\omega\neq0$.
        \begin{itemize}
            \item[(a)] If $g(\psi,\psi)\neq0$ on $W$, then $M$ is locally isomorphic, as a para-$C_{12}$-manifold, to the product of a three-dimensional $\eta$-Einstein para-$C_{12}$-manifold $(\Sigma\times I,\,h_\Sigma+a^2dt^2)$, with $\Sigma$ a Lorentzian surface and $a=a(\Sigma)$, and a para-K\"ahler--Einstein manifold $(K,h_K)$ with $\Ric_K=\lambda h_K$. The same holds for $\tilde g$.
            \item[(b)] If $g(\psi,\psi)=0$ on $W$, then $\lambda=0$ and $\dv\psi=0$ on $W$; hence $S=0$ and $\tilde S=0$, i.e.\ $g$ and $\tilde g$ are Ricci-flat. On a dense open subset of $W$, either $\Hess_ha=0$ and $a\psi$ is a parallel null vector field of $h$, or the components $\psi_\pm$ of $\psi$ in $T^\pm$ are both non-zero and $\Hess_ha=c\,a^2\,\psi_+^\flat\sym\psi_-^\flat$. In the latter case, whenever $M$ is locally a product of a para-$C_{12}$-manifold $M_1$, with $\xi$ tangent to $M_1$, and a para-K\"ahler manifold, then $\dim M_1\geq5$; in particular no three-dimensional factor as in (a) splits off.
            \item[(c)] The union of the interiors of $\{g(\psi,\psi)\neq0\}$ and $\{g(\psi,\psi)=0\}$ is dense in $W$.
        \end{itemize}
    \end{corollary}

    \begin{proof}
        As in Corollary~\ref{cor:R}, $a$ satisfies \eqref{eq:QE} with $m=1$ on the para-K\"ahler base, and $g(\psi,\psi)=a^{-2}h(\nabla a,\nabla a)$. Part (a) follows from Theorem~\ref{thm:QE}(ii) exactly as in Corollary~\ref{cor:R}. In (b), Theorem~\ref{thm:QE}(iii) gives $\lambda=0$ and $\Delta_ha=0$; by Lemma~\ref{lem:Ricmodel}, $\dv\psi=\Delta_ha/a=0$, so $S=0$ by \eqref{eq:etaE}, and $\tilde S=S+2\dv\psi\,\eta\otimes\eta=0$ by Theorem~\ref{thm:compare}(ii). The description of $\Hess_ha$ is Theorem~\ref{thm:QE}(iii-1)/(iii-2) with $\zeta=\nabla a=a\psi$. If $M=M_1\times K'$ locally with $\xi$ tangent to $M_1$, then $M_1$ is totally geodesic, so $\psi=-\nabla_\xi\xi$ is tangent to $M_1$, and the leaves of $D$ split as $N_1\times K'$ with $TN_1$ non-degenerate, $P$-invariant and containing $\zeta$; by Theorem~\ref{thm:QE}(iii-2), $\dim N_1\geq4$, i.e.\ $\dim M_1\geq5$. Part (c) holds because, for every continuous function $F$ on $W$, the set $\{F\neq0\}\cup\operatorname{int}\{F=0\}$ is open and dense in $W$.
    \end{proof}

    \begin{corollary}[Classification in setting (R)]\label{cor:classR}
        Let $(M^{2n+1},\varphi,\xi,\eta,g)$ be an $\eta$-Einstein $C_{12}$-manifold with constant $\lambda$. Every point of a dense open subset of $M$ has a neighbourhood which, as a $C_{12}$-manifold, is of one of the following types:
        \begin{itemize}
            \item[(a)] normal: $g=h+dt^2$ with $\Ric_h=\lambda h$;
            \item[(b)] the product of a three-dimensional $\eta$-Einstein $C_{12}$-manifold $(\Sigma\times I,h_\Sigma+a^2dt^2)$, $a=a(\Sigma)$, $da\neq0$, and a K\"ahler--Einstein manifold $(K,h_K)$ with $\Ric_K=\lambda h_K$;
            \item[(c)] $\lambda\neq0$, and the product of a three-dimensional space form of curvature $\lambda/2$, carrying a $C_{12}$-structure with non-closed $\omega$, and a K\"ahler--Einstein manifold with $\Ric_K=\lambda h_K$;
            \item[(d)] $\lambda=0$, and the product of a flat $C_{12}$-manifold with non-closed $\omega$ and a Ricci-flat K\"ahler manifold.
        \end{itemize}
        In (b)--(d) the K\"ahler factor is absent if $n=1$. If $M$ is connected and a point of type (c) or (d) exists, then $S=\lambda g$ on $M$.
    \end{corollary}

    \begin{proof}
        The sets $\operatorname{int}\{\omega=0\}$, $\{\omega\neq0\}\cap\operatorname{int}\{d\omega=0\}$ and $\{d\omega\neq0\}$ are open, and their union is dense in $M$. On the first, the structure is normal, $\sigma=1$ after reparametrising $t$ (Theorem~\ref{thm:twisted}), and \eqref{eq:base} reads $\Ric_h=\lambda h$: type (a). On the second, Corollary~\ref{cor:R} gives type (b); for $n=1$ this is Theorem~\ref{thm:twisted} with $\sigma=a(x)$. On a dense open subset of the third, Theorem~\ref{thm:nonclosedR}(a) gives type (c) or (d). The last assertion is Theorem~\ref{thm:nonclosed}(iv).
    \end{proof}

    \begin{remark}\label{rem:RvsP}
        In setting (P) the splitting of Corollary~\ref{cor:R} persists where $\psi$ is not null. Where $\psi$ is null and non-zero and $\omega$ is closed, an $\eta$-Einstein para-$C_{12}$-manifold is automatically Ricci-flat, and it need not split, even when it is not flat (Example~\ref{ex:iso5}); in dimension five such manifolds with $\psi_\pm\neq0$ are described by \cite[Theorem~4.1(a)]{SS-QE} with $m=1$ and $u=a$, and in dimension $2n+1\geq7$ by \cite[Theorem~5.1(a)]{SS-QE}, in closed form for $2n+1=7$ by \cite[Corollary~5.4(a)]{SS-QE}. Thus, in setting (P), Corollary~\ref{cor:P} and Theorem~\ref{thm:nonclosedR}(b) give a partial classification; the open cases are listed in Section~\ref{sec:open}.
    \end{remark}
    \section{Ricci solitons of the $\eta$-reversed metric}\label{sec:solitons}
    \noindent
    For a pseudo-Riemannian metric $\hat g$ with Ricci tensor $\hat S$, a vector field $V$ and $\lambda\in\R$ put $\mathcal E_{V,\lambda}[\hat g]:=\Lie_V\hat g+2\hat S-2\lambda\hat g$; thus $(\hat g,V,\lambda)$ is a Ricci soliton if and only if $\mathcal E_{V,\lambda}[\hat g]=0$. For fixed $\lambda$ we abbreviate $\mathcal E_V:=\mathcal E_{V,\lambda}[g]$ and $\tilde{\mathcal E}_V:=\mathcal E_{V,\lambda}[\tilde g]$. All results of this section hold in setting (Tw) of Remark~\ref{rem:settingE}; the almost (para)contact structure is used only to describe the base in Corollary~\ref{cor:Vxi}.

    \begin{lemma}\label{lem:Lie}
        \begin{itemize}
            \item[(a)] For $V=f\xi$: $(\Lie_{f\xi}\tilde g)(X,Y)=-X(f)\eta(Y)-Y(f)\eta(X)+f\bigl\{\eta(X)\omega(Y)+\eta(Y)\omega(X)\bigr\}$.
            \item[(b)] If $\eta(W)=0$, then $\Lie_W\eta=\omega(W)\eta$ and $\Lie_W\tilde g=\Lie_Wg-4\omega(W)\,\eta\otimes\eta$.
        \end{itemize}
    \end{lemma}

    \begin{proof}
        (a) $\nt_X(f\xi)=X(f)\xi+f\eta(X)\psi$ by Theorem~\ref{thm:compare}(iv), and $\tilde g(\xi,\cdot)=-\eta$, $\tilde g(\psi,\cdot)=\omega$.
        (b) $(\Lie_W\eta)X=(\nabla_W\eta)X+\eta(\nabla_XW)=0+g(\nabla_XW,\xi)=-g(W,\nabla_X\xi)=\eta(X)\omega(W)$; then $\Lie_W\tilde g=\Lie_W(g-2\eta\otimes\eta)$.
    \end{proof}

    \begin{theorem}[General potential]\label{thm:general}
        Let $V=f\xi+W$ with $\eta(W)=0$. Then $(\tilde g,V,\lambda)$ is a Ricci soliton if and only if
        \begin{itemize}
            \item[(a)] $(\Lie_Wg+2S-2\lambda g)|_{D\times D}=0$,
            \item[(b)] $X(f)=f\omega(X)+g(\nabla_\xi W,X)$ for all $X\in D$,
            \item[(c)] $\xi(f)=\lambda+\dv\psi-\omega(W)$.
        \end{itemize}
    \end{theorem}

    \begin{proof}
        By Lemma~\ref{lem:Lie}, Theorem~\ref{thm:compare}(ii) and $\tilde g=g-2\eta\otimes\eta$,
        \begin{equation}\label{eq:Tid}
            \tilde{\mathcal E}_V=\mathcal E_W-df\otimes\eta-\eta\otimes df+f(\eta\otimes\omega+\omega\otimes\eta)+4\bigl(\lambda+\dv\psi-\omega(W)\bigr)\eta\otimes\eta .
        \end{equation}
        On $D\times D$ this is $\mathcal E_W$, which gives (a). For $X\in D$ we have $\mathcal E_W(X,\xi)=g(\nabla_XW,\xi)+g(\nabla_\xi W,X)=g(\nabla_\xi W,X)$, as $g(\nabla_XW,\xi)=\eta(X)\omega(W)=0$; hence $\tilde{\mathcal E}_V(X,\xi)=g(\nabla_\xi W,X)-X(f)+f\omega(X)$, which gives (b). Finally $\mathcal E_W(\xi,\xi)=2\omega(W)-2\dv\psi-2\lambda$ and $\tilde{\mathcal E}_V(\xi,\xi)=2(\lambda+\dv\psi-\omega(W))-2\xi(f)$, which gives (c).
    \end{proof}

    \begin{corollary}[$V\perp\xi$]\label{cor:Vperp}
        If $\eta(V)=0$, then $(\tilde g,V,\lambda)$ is a Ricci soliton if and only if $(g,V,\lambda)$ is, and then $\lambda=\omega(V)-\dv\psi$; in particular this function is constant.
    \end{corollary}

    \begin{proof}
        Take $f=0$ and $W=V$ in \eqref{eq:Tid}: $\tilde{\mathcal E}_V=\mathcal E_V+4(\lambda+\dv\psi-\omega(V))\eta\otimes\eta$ and $\mathcal E_V(\xi,\xi)=2\omega(V)-2\dv\psi-2\lambda$. If $\tilde{\mathcal E}_V=0$, then $\tilde{\mathcal E}_V(\xi,\xi)=2(\lambda+\dv\psi-\omega(V))=0$, whence $\tilde{\mathcal E}_V=\mathcal E_V=0$; conversely, $\mathcal E_V=0$ gives $\mathcal E_V(\xi,\xi)=0$ and again $\tilde{\mathcal E}_V=\mathcal E_V$.
    \end{proof}

    \noindent
    Corollary~\ref{cor:Vperp} produces solitons of $\tilde g$ with $W\neq0$ that are neither trivial nor products with $\omega=0$; see Example~\ref{ex:cigar}, which is built on Hamilton's cigar.

    \begin{theorem}[$V=f\xi$]\label{thm:fxi}
        $(\tilde g,f\xi,\lambda)$ is a Ricci soliton if and only if
        \begin{itemize}
            \item[(i)] $S=\lambda g-(\lambda+\dv\psi)\,\eta\otimes\eta$, i.e.\ $M$ is $\eta$-Einstein with constant $\lambda$, and
            \item[(ii)] $df=f\omega+(\lambda+\dv\psi)\,\eta$.
        \end{itemize}
    \end{theorem}

    \begin{proof}
        This is Theorem~\ref{thm:general} with $W=0$: (i) is (a) together with $S(\xi,\cdot)=-\dv\psi\,\eta$ (Lemma~\ref{lem:basic}(iv)), and (ii) is (b) together with (c).
    \end{proof}

    \begin{theorem}[Ricci solitons with potential $f\xi$]\label{thm:rigidity}
        Let $M$ be connected, let $f$ be a smooth function and $\lambda\in\R$, and suppose that $(\tilde g,f\xi,\lambda)$ is a Ricci soliton. Then exactly one of the following holds.
        \begin{itemize}
            \item[(A)] $\lambda\neq0$ and $\omega=0$ on $M$, i.e.\ $\xi$ is parallel; $S=\lambda(g-\eta\otimes\eta)$, and locally $\tilde g=-d\tau^2+h$ with $\Ric_h=\lambda h$ and $f\xi=(\lambda\tau+k_0)\partial_\tau$, $k_0\in\R$. In particular $f\xi$ is not a Killing field of $\tilde g$.
            \item[(B)] $S=\lambda g$ and $\tilde S=\lambda\tilde g$ on $M$, and $f\xi$ is a Killing field of $\tilde g$. If $f\not\equiv0$, then $f$ has no zeros and $\omega=d\ln|f|$; in particular $\omega$ is exact, and locally $g=h+a(x)^2d\tau^2$ is a warped product with $f\xi=k_0\partial_\tau$.
        \end{itemize}
    \end{theorem}

    \begin{proof}
        By Theorem~\ref{thm:fxi}, $M$ is $\eta$-Einstein with constant $\lambda$, $S=\lambda g-\Theta\,\eta\otimes\eta$ with $\Theta=\lambda+\dv\psi$, and $df=f\omega+\Theta\eta$. By Theorem~\ref{thm:nonclosed}(iii), either $\Theta\equiv0$, or $\Theta$ has no zeros and $\omega$ is closed.

        Suppose that $\Theta$ has no zeros. Applying $d$ to $df=f\omega+\Theta\eta$ and using $d\omega=0$, $d\eta=\omega\wedge\eta$ and $d\Theta=-2\Theta\omega$ (Theorem~\ref{thm:nonclosed}(ii)), we obtain
        \[
            0=df\wedge\omega+d\Theta\wedge\eta+\Theta\,d\eta=\Theta\,\eta\wedge\omega-2\Theta\,\omega\wedge\eta+\Theta\,\omega\wedge\eta=2\Theta\,\eta\wedge\omega .
        \]
        Hence $\eta\wedge\omega=0$, and since $\omega(\xi)=0$, $\omega=0$. Then $\psi=0$, $\Theta=\lambda\neq0$ and $S=\lambda(g-\eta\otimes\eta)$. Locally $\sigma=1$ after reparametrising $t=:\tau$ (Theorem~\ref{thm:twisted}), so $\eta=d\tau$, $\xi=\partial_\tau$, $\tilde g=-d\tau^2+h$, $\Ric_h=S|_{D\times D}=\lambda h$ by Lemma~\ref{lem:Ricmodel}, and $df=\lambda\,d\tau$ gives $f=\lambda\tau+k_0$. By Lemma~\ref{lem:Lie}(a), $\Lie_{f\xi}\tilde g=-2\lambda\,d\tau^2\neq0$. This is (A).

        Suppose that $\Theta\equiv0$. Then $S=\lambda g$, $\dv\psi=-\lambda$ and $\tilde S=S+2\dv\psi\,\eta\otimes\eta=\lambda\tilde g$ by Theorem~\ref{thm:compare}(ii). Moreover $df=f\omega$, so Lemma~\ref{lem:Lie}(a) gives $\Lie_{f\xi}\tilde g=-(df-f\omega)\otimes\eta-\eta\otimes(df-f\omega)=0$. If $f\not\equiv0$, then, as in the proof of Theorem~\ref{thm:nonclosed}(iii), the linear equation $df=f\omega$ shows that $f$ has no zeros, and $\omega=d\ln|f|$. By Theorem~\ref{thm:twisted} (closedness), locally $\sigma=a(x)$ after reparametrising $t=:\tau$, and $d\ln|f|=\omega=d_N\ln a$ gives $f=k_0a$, i.e.\ $f\xi=k_0\partial_\tau$. This is (B). The two cases are exclusive because $\Theta=\lambda\neq0$ in (A).
    \end{proof}

    \begin{remark}\label{rem:sharp}
        (1) For $f\equiv0$, i.e.\ when $\tilde g$ is Einstein, $\omega$ need not be closed: see Example~\ref{ex:timedep}.
        (2) Theorem~\ref{thm:rigidity} contains the local statement that, on every connected open set where $f\neq0$, the twisting function $\sigma$ of the model \eqref{eq:model} separates variables, i.e.\ the twisted product is locally a warped product; in case (B) one has $f=k_0a$, and in case (A) $a=1$.
        (3) Ricci solitons on warped products $B\times_uF$ have been studied mainly for gradient potentials, both in the Riemannian and in the semi-Riemannian setting \cite{BT22,FFG17}. The potentials $f\xi$ above are in general not gradient fields: by Theorem~\ref{thm:fxi}(ii), $d\bigl(\tilde g(f\xi,\cdot)\bigr)=-d(f\eta)=-df\wedge\eta-f\,d\eta=-2f\,\omega\wedge\eta$, which vanishes only where $f\omega=0$. Moreover, our base is a twisted product with a one-dimensional fibre, and Theorem~\ref{thm:rigidity} shows that the twisting is forced to be a warping.
    \end{remark}

    \begin{corollary}[$V=\xi$]\label{cor:Vxi}
        $(\tilde g,\xi,\lambda)$ is a Ricci soliton if and only if $\omega=0$, $\lambda=0$ and $\Ric_h=0$; then locally $\tilde g=-dt^2+h$ with $h$ Ricci-flat and K\"ahler (resp.\ para-K\"ahler).
    \end{corollary}

    \begin{proof}
        Theorem~\ref{thm:fxi}(ii) with $f=1$ gives $0=\omega+(\lambda+\dv\psi)\eta$, so $\omega=0$; hence $\psi=0$, $\dv\psi=0$ and $\lambda=0$. By Theorem~\ref{thm:twisted} we may take $\sigma=1$, and Theorem~\ref{thm:fxi}(i) gives $\Ric_h=0$.
    \end{proof}
    \section{Examples}\label{sec:examples}
    \subsection{Setting (R)}
    \noindent
    \begin{example}\label{ex:S2}
        Let $a\equiv1$ and let $h$ be the round unit sphere $S^2$, so that $\Ric_h=h$. Then $\tilde g=-dt^2+h$, and $V=(t+k_0)\partial_t$ is a Ricci soliton vector field with $\lambda=1$ (Theorem~\ref{thm:fxi}, $\dv\psi=0$). This is case (A) of Theorem~\ref{thm:rigidity}.
    \end{example}

    \begin{example}[$S^3$ and de Sitter space]\label{ex:S3dS}
        Let $a=\cos r$ and $h=dr^2+\sin^2r\,d\theta^2$ on $0<r<\pi/2$. Then $g=\cos^2r\,dt^2+dr^2+\sin^2r\,d\theta^2$ is the round metric of $S^3$, so $S=2g$, and it carries a non-normal $C_{12}$-structure with $\omega=d\ln\cos r\neq0$; $\tilde g$ is the static patch of three-dimensional de Sitter space, $\tilde S=2\tilde g$. Here $\dv\psi=-2$ and $\lambda=2$. Since $\partial_t$ is a Killing field of $\tilde g$, $(\tilde g,k_0\partial_t,2)$ is a Ricci soliton with $f=k_0\cos r$; this is case (B) of Theorem~\ref{thm:rigidity}.
    \end{example}

    \begin{example}[Time-dependent lapse]\label{ex:timedep}
        Let $h$ be the round unit sphere $S^2\subset\R^3$, with position vector $X$, and let $\sigma(X,t)=\langle\mathbf{q}(t),X\rangle$ for a smooth curve $\mathbf{q}$ in $\R^3$, on the open set where $\sigma>0$. Since $\Hess_h\sigma=-\sigma h$, Lemma~\ref{lem:Ricmodel} gives $S=2h$ on $D$ and $S(\xi,\xi)=-\Delta_h\sigma/\sigma=2$, so $S=2g$, and Theorem~\ref{thm:compare}(ii) gives $\tilde S=2g-4\eta\otimes\eta=2\tilde g$. On the other hand, $d\omega=dt\wedge d_N\bigl(\langle\mathbf{q}',X\rangle/\langle\mathbf{q},X\rangle\bigr)$ vanishes if and only if $\mathbf{q}'\parallel\mathbf{q}$. This is type (a1) of Theorem~\ref{thm:nonclosedR} with $\lambda=2$; accordingly, $\Theta=\scal_h-\lambda=0$, as required by Theorem~\ref{thm:nonclosed}. Since $(\tilde g,0,2)$ is a Ricci soliton, the example also shows that the exactness of $\omega$ in Theorem~\ref{thm:rigidity}(B) needs $f\not\equiv0$.
    \end{example}

    \begin{example}[Non-closed $\omega$, flat]\label{ex:flatnc}
        Let $\R^{2n}$ carry the flat K\"ahler structure with $P\partial_{x^{2j-1}}=\partial_{x^{2j}}$, let $1\leq k\leq n$, and let $\sigma=1+\sum_{j=1}^kt^j\,x^{2j-1}$ on the set where $\sigma>0$; for $n=k=1$ this is $g=dx^2+dy^2+(1+tx)^2dt^2$. Each $\sigma_t$ is affine, so $\Hess_h\sigma_t=0$ and \eqref{eq:base} holds with $\lambda=0$. Since $\partial_t\ln\sigma=\sum_jj\,t^{j-1}x^{2j-1}/\sigma$, the form $\omega$ is not closed. The gradients $\nabla\sigma_t=\sum_jt^j\partial_{x^{2j-1}}$, $t\in I$, span $\Span\{\partial_{x^1},\partial_{x^3},\dots,\partial_{x^{2k-1}}\}$ (a Vandermonde argument), so the distribution $E_0$ in the proof of Theorem~\ref{thm:nonclosedR} has rank $2k$: this is type (a2) with the given $k$ and $K=\R^{2n-2k}$, and $g$ is flat. By Theorem~\ref{thm:compare}(v), the $\eta$-reversed metric $\tilde g$ is flat as well.
    \end{example}

    \begin{example}[Non-closed $\omega$, hyperbolic]\label{ex:hypnc}
        Let $h=dr^2+\sinh^2r\,d\theta^2$ be the hyperbolic plane and $\sigma=\cosh r+t\sinh r\cos\theta$, $|t|<1$. As $\sigma(\cdot,t)$ is the restriction of a linear function to the hyperboloid model, $\Hess_h\sigma=\sigma h$, so \eqref{eq:base} holds with $\lambda=-2$ and $d\omega\neq0$. By Theorem~\ref{thm:nonclosed}, $S=-2g$: $(H^2\times I,h+\sigma^2dt^2)$ is locally real hyperbolic space with a $C_{12}$-structure for which $\omega$ is not closed (type (a1) with $\lambda=-2$).
    \end{example}

    \begin{example}[Product]\label{ex:product}
        Let $\Sigma=S^2$ with $a=\cos r$ as in Example~\ref{ex:S3dS}, and let $K=S^2(1/\sqrt2)$. Then
        \[
            g=\cos^2r\,dt^2+dr^2+\sin^2r\,d\theta^2+\tfrac12\bigl(ds^2+\sin^2s\,dw^2\bigr)
        \]
        is $S^3\times S^2(1/\sqrt2)$ with $S=2g$, and $\tilde g$ is $\mathrm{dS}_3\times S^2(1/\sqrt2)$, again Einstein. This is the shape predicted by Corollary~\ref{cor:R}.
    \end{example}

    \begin{example}[A soliton with potential orthogonal to $\xi$]\label{ex:cigar}
        Over the flat K\"ahler half-plane $\{(x,y):x>0\}$ let $\sigma=\tanh x$, i.e.\ $g=dx^2+dy^2+\tanh^2x\,dt^2$. Then $\psi=(2/\sinh2x)\,\partial_x\neq0$ and $\omega=d\ln\tanh x$ is closed. For $t\in\R/2\pi\mathbb{Z}$, $dx^2+\tanh^2x\,dt^2$ is Hamilton's cigar \cite{Ham88} in geodesic polar coordinates, so $(M,g)$ is the product of the cigar without its tip and a line. The function $F=-2\ln\cosh x$ satisfies $\Ric_g+\Hess_gF=0$, i.e.\ $(g,\nabla F,0)$ is a steady gradient Ricci soliton, and $V:=\nabla F=-2\tanh x\,\partial_x$ is orthogonal to $\xi$. By Corollary~\ref{cor:Vperp}, $(\tilde g,V,0)$ is a Ricci soliton of $\tilde g=dx^2+dy^2-\tanh^2x\,dt^2$; accordingly $\omega(V)-\dv\psi=-2\cosh^{-2}x+2\cosh^{-2}x=0$. Here $V$ is not a Killing field of $\tilde g$, $W=V\neq0$, and $g$ is not Einstein.
    \end{example}

    \subsection{Setting (P)}
    \noindent
    In the following examples $M=N\times I$ carries the para-$C_{12}$-structure of Theorem~\ref{thm:twisted} over the given para-K\"ahler base.

    \begin{example}[Non-isotropic, Einstein]\label{ex:dSP}
        Let $\Sigma$ carry $h_\Sigma=-d\tau^2+\cosh^2\tau\,d\theta^2$ with $P_\Sigma\partial_\tau=(\cosh\tau)^{-1}\partial_\theta$, and let $a=\sinh\tau$, $\tau>0$. Then $\Hess_\Sigma a=-a\,h_\Sigma$ and $h_\Sigma(\nabla a,\nabla a)=-\cosh^2\tau\neq0$. Let $K$ carry $h_K=\frac12(-ds^2+\cosh^2s\,dw^2)$ with $P_K\partial_s=(\cosh s)^{-1}\partial_w$, so that $\Ric_K=2h_K$. Then $N=\Sigma\times K$ satisfies \eqref{eq:base} with $\lambda=2$, and $g=-d\tau^2+\cosh^2\tau\,d\theta^2+\sinh^2\tau\,dt^2+h_K$ is an Einstein para-$C_{12}$-metric with $S=2g$ and $\tilde S=2\tilde g$: locally $\mathrm{dS}_3\times\mathrm{dS}_2$, the paracontact analogue of Example~\ref{ex:product}.
    \end{example}

    \begin{example}[Non-isotropic, $\eta$-Einstein but not Einstein]\label{ex:ODE}
        Let $a(r)>0$, $a'(r)>0$ be a solution of $a\,a''=\frac\lambda2a^2+C_1$ (for instance $\lambda=0$, $C_1=1$, $a(0)=a'(0)=1$, near $r=0$), and let $h_\Sigma=-dr^2+a'(r)^2d\theta^2$ with $P_\Sigma\partial_r=a'(r)^{-1}\partial_\theta$. Then $\Hess_\Sigma a=-a''h_\Sigma$, $\Ric_\Sigma-\Hess_\Sigma a/a=\lambda h_\Sigma$ and $\mu=-2C_1$. Hence $(\Sigma\times I,\,h_\Sigma+a^2dt^2)$ is an $\eta$-Einstein para-$C_{12}$-manifold with $S=\lambda g+2C_1a^{-2}\,\eta\otimes\eta$, which is not Einstein if $C_1\neq0$; in accordance with Proposition~\ref{prop:constcoef}, the coefficient of $\eta\otimes\eta$ is then not constant. Products with para-K\"ahler--Einstein manifolds $K$, $\Ric_K=\lambda h_K$, give examples of type (a) of Corollary~\ref{cor:P} in every dimension.
    \end{example}

    \begin{example}[Non-closed $\omega$, de Sitter]\label{ex:dSnc}
        Let $h_\Sigma=-d\tau^2+\cosh^2\tau\,d\theta^2$ with $P_\Sigma$ as in Example~\ref{ex:dSP}, and $\sigma=\sinh\tau+t\cosh\tau\cos\theta$ on the set where $\sigma>0$. As $\sigma(\cdot,t)$ is the restriction of a linear function to the hyperboloid model of the de Sitter plane, $\Hess_{h_\Sigma}\sigma=-\sigma h_\Sigma$; hence \eqref{eq:base} holds with $\lambda=2$, $d\omega\neq0$, and $S=2g$. This is the first case of Theorem~\ref{thm:nonclosedR}(b): $(\Sigma\times I,h_\Sigma+\sigma^2dt^2)$ is locally three-dimensional de Sitter space, with a para-$C_{12}$-structure for which $\omega$ is not closed.
    \end{example}

    \begin{example}[Non-closed $\omega$, Ricci-flat, not flat]\label{ex:ncP}
        Let $h=dx^1\sym dy^1+dx^2\sym dy^2+2x^2y^1\,dx^2\sym dy^1$ be the para-K\"ahler metric on $\R^4$ of \cite[Example~6.3]{SS-QE}, with $P\partial_{x^i}=\partial_{x^i}$ and $P\partial_{y^i}=-\partial_{y^i}$, and let $\sigma=1+t\,y^1$ on the set where $\sigma>0$. Since $\nabla y^1=\partial_{x^1}$ is parallel and $\Ric_h=0$, we have $\Hess_h\sigma_t=0$, so \eqref{eq:base} holds with $\lambda=0$, and $M$ is an $\eta$-Einstein para-$C_{12}$-manifold with $S=0$. As $\partial_t\ln\sigma=y^1/(1+ty^1)$, the form $\omega$ is not closed, and all gradients $\nabla\sigma_t=t\,\partial_{x^1}$ are null: this is the last case of Theorem~\ref{thm:nonclosedR}(b). The metric $g$ is Ricci-flat but not flat, because its totally geodesic leaves $N\times\{t\}$ are isometric to the non-flat manifold $(N,h)$.
    \end{example}

    \begin{example}[Isotropic, five-dimensional]\label{ex:iso5}
        Let $(h,u)$ be one of the following para-K\"ahler metrics and functions on open subsets of $\R^4$, with $P\partial_{x^i}=\partial_{x^i}$ and $P\partial_{y^i}=-\partial_{y^i}$ \cite[Examples~6.1--6.3]{SS-QE}: (i) $h=dx^1\sym dy^1+dx^2\sym dy^2$ and $u=x^1+y^2$; (ii) $h=2v\bigl(dx^1\sym dy^1+dx^2\sym dy^2\bigr)+2(x^2+y^1)\,dx^1\sym dy^2$ and $u=v^2$, where $v=x^1+y^2>0$; (iii) $h=dx^1\sym dy^1+dx^2\sym dy^2+2x^2y^1\,dx^2\sym dy^1$ and $u=y^1>0$. Let $g=h+a^2dt^2$ with $a:=u$. Since \eqref{eq:QE} holds with $m=1$ and $\lambda=0$, $M$ is an $\eta$-Einstein para-$C_{12}$-manifold with closed $\omega=d\ln a\neq0$ and $g(\psi,\psi)=a^{-2}h(\nabla a,\nabla a)=0$; by Corollary~\ref{cor:P}(b), $g$ and $\tilde g$ are Ricci-flat. Case (i) gives a flat manifold and case (ii) a non-flat one; in both cases $\psi_\pm\neq0$, so by Corollary~\ref{cor:P}(b) no three-dimensional para-$C_{12}$-factor splits off, in contrast with Corollary~\ref{cor:R}. Case (iii) gives a Ricci-flat, non-flat para-$C_{12}$-manifold with $\psi\in T^+$.
    \end{example}
    \section{Left-invariant structures on three-dimensional Lie groups}\label{sec:lie}
    \noindent
    We use the three-dimensional non-abelian real Lie algebras of Patera et al.\ \cite{PSWZ76}, with the non-zero brackets
    \begin{align*}
        &A_{3,1}\colon [e_2,e_3]=e_1;\qquad A_{3,2}\colon [e_1,e_3]=e_1,\ [e_2,e_3]=e_1+e_2;\\
        &A_{3,3}\colon [e_1,e_3]=e_1,\ [e_2,e_3]=e_2;\qquad A_{3,4}\colon [e_1,e_3]=e_1,\ [e_2,e_3]=-e_2;\\
        &A^{\gamma}_{3,5}\colon [e_1,e_3]=e_1,\ [e_2,e_3]=\gamma e_2\ (0<|\gamma|<1);\\
        &A_{3,6}\colon [e_1,e_3]=-e_2,\ [e_2,e_3]=e_1;\\
        &A^{\gamma}_{3,7}\colon [e_1,e_3]=\gamma e_1-e_2,\ [e_2,e_3]=e_1+\gamma e_2\ (\gamma>0);\\
        &A_{3,8}\colon [e_1,e_3]=-2e_2,\ [e_1,e_2]=e_1,\ [e_2,e_3]=e_3;\\
        &A_{3,9}\colon [e_1,e_2]=e_3,\ [e_2,e_3]=e_1,\ [e_3,e_1]=e_2.
    \end{align*}
    Recall that $A_{3,1}$ is the Heisenberg algebra, $A_{3,4}\cong\mathfrak e(1,1)$, $A_{3,6}\cong\mathfrak e(2)$, $A_{3,8}\cong\mathfrak{sl}(2,\R)$ and $A_{3,9}\cong\mathfrak{su}(2)$ are unimodular, and $A_{3,2}$, $A_{3,3}$, $A^\gamma_{3,5}$, $A^\gamma_{3,7}$ are not. Let $g$ be an \emph{arbitrary} left-invariant Riemannian metric on the simply connected Lie group with one of these Lie algebras $\mathfrak g$. In dimension three, $\varphi$ is necessarily the rotation by $\pm\pi/2$ in $\xi^\perp$, and by Theorem~\ref{thm:twisted} and Corollary~\ref{cor:dim3} a left-invariant almost contact metric structure is $C_{12}$ if and only if $\nabla_X\xi=-\eta(X)\psi$ (cf.\ \cite[Thm.~3.3]{Bel24}). The following lemma turns this into an algebraic condition.

    \begin{lemma}\label{lem:liecrit}
        Let $\xi\in\mathfrak g$ be a unit vector and $D=\xi^\perp$. Then $\nabla_X\xi=-\eta(X)\psi$ for all $X$ if and only if (a) $D$ is a subalgebra of $\mathfrak g$ and (b) $g([\xi,X],X)=0$ for all $X\in D$. In this case $\omega(Z)=g([\xi,Z],\xi)$ and $\dv\psi=-\tr\operatorname{ad}_\psi$.
    \end{lemma}

    \begin{proof}
        For left-invariant vector fields, $2g(\nabla_XY,Z)=g([X,Y],Z)-g([Y,Z],X)+g([Z,X],Y)$. For $X,Y\in D$ this gives
        \begin{align*}
            2g(\nabla_X\xi,Y)+2g(\nabla_Y\xi,X)&=-2g([\xi,X],Y)-2g([\xi,Y],X),\\
            2g(\nabla_X\xi,Y)-2g(\nabla_Y\xi,X)&=-2\eta([X,Y]).
        \end{align*}
        Since $g(\nabla_X\xi,\xi)=0$, the condition $\nabla_X\xi=0$ for all $X\in D$, which is equivalent to $\nabla_X\xi=-\eta(X)\psi$ for all $X$, holds if and only if both expressions vanish, i.e.\ if and only if (a) and (b) hold. Moreover $g(\nabla_\xi\xi,Z)=g([Z,\xi],\xi)$, so $\omega(Z)=g([\xi,Z],\xi)$, and, for a $g$-orthonormal basis $(e_i')$, $\dv X=\sum_ig(\nabla_{e_i'}X,e_i')=\sum_ig([e_i',X],e_i')=-\tr\operatorname{ad}_X$ for every left-invariant $X$.
    \end{proof}

    \begin{proposition}\label{prop:lieclass}
        Let $g$ be an arbitrary left-invariant metric.
        \begin{itemize}
            \item[(a)] Let $\mathfrak g$ be unimodular, and let $(f_1,f_2,f_3)$ be a $g$-orthonormal Milnor frame \cite{Mil76}, i.e.\ $[f_2,f_3]=l_1f_1$, $[f_3,f_1]=l_2f_2$, $[f_1,f_2]=l_3f_3$, numbered so that $l_1\geq l_2\geq l_3$ and oriented so that at most one $l_i$ is negative (reversing the orientation changes the signs of all $l_i$, and the condition below is invariant under this change); equivalently, $[u,v]=\Lambda(u\times v)$ with $\Lambda=\operatorname{diag}(l_1,l_2,l_3)$ and the cross product $\times$ of the oriented orthonormal frame $(f_i)$. A unit vector $\xi\in\mathfrak g$ defines a left-invariant $C_{12}$-structure if and only if
            \[
                l_2=l_1+l_3\qquad\text{and}\qquad \xi=\epsilon_1\sqrt{\frac{-l_3}{l_1-l_3}}\,f_1+\epsilon_3\sqrt{\frac{l_1}{l_1-l_3}}\,f_3,\quad \epsilon_1,\epsilon_3\in\{\pm1\}.
            \]
            Then $\psi=(\Lambda\xi)\times\xi=-\epsilon_1\epsilon_3\sqrt{-l_1l_3}\,f_2$ and $\dv\psi=0$. Consequently there are no such structures on $A_{3,1}$ and $A_{3,9}$; on $A_{3,6}$ they exist if and only if $l_1=l_2$, i.e.\ $g$ is flat, and then $\xi=\pm f_3$ and $\psi=0$; on $A_{3,4}$ they exist if and only if $l_1=-l_3$ (and $l_2=0$); on $A_{3,8}$ they exist if and only if $l_2=l_1+l_3$.
            \item[(b)] Let $\mathfrak g$ be one of the non-unimodular algebras $A_{3,2}$, $A_{3,3}$, $A^\gamma_{3,5}$, $A^\gamma_{3,7}$, write $\mathfrak g=\mathfrak n\rtimes\R e_3$ with $\mathfrak n=\Span\{e_1,e_2\}$, and put $L:=\operatorname{ad}_{e_3}|_{\mathfrak n}$. Left-invariant $C_{12}$-structures exist on $A_{3,3}$ for every unit vector $\xi\in\mathfrak n$, and on $A^\gamma_{3,5}$ if and only if $g(e_1,e_2)=0$, in which case $\xi\in\{\pm e_1/|e_1|,\pm e_2/|e_2|\}$; there are none on $A_{3,2}$ and $A^\gamma_{3,7}$. In all cases $L\xi=\nu\xi$ for some $\nu\in\R\setminus\{0\}$, and
            \[
                \psi=-\nu\,(e^3)^\sharp,\qquad \dv\psi=\nu\,|e^3|^2\tr L ,
            \]
            where $e^3\in\mathfrak g^*$ is given by $e^3(e_3)=1$, $e^3(\mathfrak n)=0$.
        \end{itemize}
    \end{proposition}

    \begin{proof}
        (a) For $X,Y\in\xi^\perp$ the vector $X\times Y$ is a multiple of $\xi$, so $\xi^\perp$ is a subalgebra if and only if $g(\Lambda\xi,\xi)=0$. The map $J:=\xi\times\cdot$ is a rotation by $\pi/2$ of $\xi^\perp$, and condition (b) of Lemma~\ref{lem:liecrit} reads $g(\Lambda JX,X)=0$ for $X\in\xi^\perp$. If $\Lambda_\xi$ denotes the compression of $\Lambda$ to $\xi^\perp$, a symmetric endomorphism of a Euclidean plane, this says that $\Lambda_\xi J$ is skew, i.e.\ that $\Lambda_\xi$ commutes with $J$, i.e.\ that $\Lambda_\xi=c\,\Id$ for some $c\in\R$. Thus, in an orthonormal basis $(\xi,X_1,X_2)$ with $\Lambda\xi\in\R X_1$, the matrix of $\Lambda$ is
        \[
            \begin{pmatrix}0&w&0\\ w&c&0\\ 0&0&c\end{pmatrix},\qquad w=|\Lambda\xi| ,
        \]
        whose eigenvalues are $c$ and the two roots of $x^2-cx-w^2$; these have sum $c$ and product $-w^2\leq0$. Hence $c=l_2=l_1+l_3$, $l_1\geq0\geq l_3$, and the $l_2$-eigenvector $X_2$ is orthogonal to $\xi$. If the $l_i$ are distinct, $X_2=\pm f_2$, and $\xi\in\Span\{f_1,f_3\}$ with $g(\Lambda\xi,\xi)=l_1\xi_1^2+l_3\xi_3^2=0$ gives the four vectors $\xi$ above. Conversely, for these vectors $\xi^\perp$ contains $f_2$, $g(\Lambda\xi,\xi)=0$, and $\Lambda_\xi$ has the eigenvalue $l_2$ on $f_2$ and trace $\tr\Lambda-g(\Lambda\xi,\xi)=2l_2$, so $\Lambda_\xi=l_2\Id$. If $l_1=l_2$, then $l_3=0$ and $g(\Lambda\xi,\xi)=l_1(\xi_1^2+\xi_2^2)=0$ forces $\xi=\pm f_3$, in accordance with the formula above; $l_2=l_3\neq0$ cannot occur, because $l_1\geq l_2$ and at most one $l_i$ is negative. Next, $\omega(Z)=g(\Lambda(\xi\times Z),\xi)=g(Z,\Lambda\xi\times\xi)$, so $\psi=\Lambda\xi\times\xi$, which is computed directly; and $\dv\psi=-\tr\operatorname{ad}_\psi=0$ by unimodularity. Finally, by Milnor's classification \cite[\S4]{Mil76}, with our orientation convention the signs of $(l_1,l_2,l_3)$ are $(+,+,+)$ for $\mathfrak{su}(2)$, $(+,0,0)$ for the Heisenberg algebra, $(+,+,0)$ for $\mathfrak e(2)$, $(+,0,-)$ for $\mathfrak e(1,1)$ and $(+,+,-)$ for $\mathfrak{sl}(2,\R)$; the condition $l_2=l_1+l_3$ is impossible in the first two cases and gives the stated conditions in the others. A metric on $\mathfrak e(2)$ is flat if and only if $l_1=l_2$ \cite{Mil76}.

        (b) Here $L$ is invertible. A two-dimensional subalgebra other than $\mathfrak n$ meets $\mathfrak n$ in a line $\R v$ and contains some $e_3+n_0$ with $n_0\in\mathfrak n$; then $[e_3+n_0,v]=Lv\in\R v$, so $Lv=\nu' v$ with $\nu'\neq0$. If $D=\mathfrak n$, write $\xi=\alpha e_3+n_1$ with $\alpha\neq0$; condition (b) of Lemma~\ref{lem:liecrit} says that $g(LX,X)=0$ for $X\in\mathfrak n$, i.e.\ that $L$ is skew for $g|_{\mathfrak n}$, which forces $\tr L=0$, impossible for non-unimodular $\mathfrak g$. Otherwise write $\xi=\alpha e_3+n_1$; (b) for $X=v$ gives $0=g([\xi,v],v)=\alpha\nu'|v|^2$, so $\alpha=0$, $\xi\in\mathfrak n$ and $\xi\perp v$. For $X=a v+b(e_3+n_0)\in D$ we get $g([\xi,X],X)=-b\,g(L\xi,X)$, so (b) holds if and only if $L\xi\perp D$, i.e.\ $L\xi\in\R\xi$. Thus $\xi$ and $v$ must be $g$-orthogonal eigenvectors of $L$ in $\mathfrak n$, and conversely every such pair gives a subalgebra $D=\xi^\perp\ni v$ satisfying (b). This is impossible for $A_{3,2}$ (one eigendirection) and for $A^\gamma_{3,7}$ (no real eigenvalues), holds for every unit $\xi\in\mathfrak n$ for $A_{3,3}$ ($L=-\Id$), and requires $g(e_1,e_2)=0$ and $\xi\in\R e_1\cup\R e_2$ for $A^\gamma_{3,5}$ ($L=-\operatorname{diag}(1,\gamma)$). By Lemma~\ref{lem:liecrit}, $\omega(e_3)=g([\xi,e_3],\xi)=-\nu$ and $\omega(\mathfrak n)=0$, so $\omega=-\nu e^3$; and $\dv\psi=\nu\tr\operatorname{ad}_{(e^3)^\sharp}=\nu\,|e^3|^2\tr L$, since $\tr\operatorname{ad}_n=0$ for $n\in\mathfrak n$.
    \end{proof}

    \begin{corollary}\label{cor:lienormal}
        Up to automorphisms of $\mathfrak g$ and homotheties, the left-invariant metrics that carry a left-invariant $C_{12}$-structure are: the metrics for which the Patera basis is orthonormal, on $A_{3,3}$, $A_{3,4}$, $A^\gamma_{3,5}$ and $A_{3,6}$; and on $A_{3,8}\cong\mathfrak{sl}(2,\R)$ the one-parameter family $g_\rho$, $0<\rho<1$, with Milnor constants $(l_1,l_2,l_3)=(1,1-\rho,-\rho)$. The metric for which the Patera basis of $A_{3,8}$ is orthonormal has Milnor constants $(2,1,-1)$, i.e.\ it is homothetic to $g_{1/2}$; the value $\rho=1$ gives the orthonormal Patera metric on $A_{3,4}$.
    \end{corollary}

    \begin{proof}
        For $A_{3,3}$, $A_{3,4}$, $A^\gamma_{3,5}$ and $A_{3,6}$ the maps $e_3\mapsto e_3+n$ ($n\in\mathfrak n$), and the linear maps of $\mathfrak n$ commuting with $L$, extended by $e_3\mapsto e_3$, are automorphisms. By Proposition~\ref{prop:lieclass} and the equivalent conditions for $A_{3,4}$ and $A_{3,6}$ obtained by applying the argument of (b) to these algebras ($g(e_1,e_2)=0$, resp.\ $g|_{\mathfrak n}$ invariant under $L$), such automorphisms bring every metric carrying a $C_{12}$-structure to the form $e^1\otimes e^1+e^2\otimes e^2+\kappa\,e^3\otimes e^3$, $\kappa>0$; the homothety by $1/\kappa$ followed by $e_i\mapsto\sqrt\kappa\,e_i$ ($i=1,2$) gives the orthonormal Patera metric. For unimodular $\mathfrak g$, two metrics with the same Milnor constants are related by the Lie algebra isomorphism that maps one Milnor frame to the other, and homotheties rescale the Milnor constants by a common positive factor; for $A_{3,8}$ we may therefore assume $l_1=1$, and then $l_2=l_1+l_3$ and $l_2>0>l_3$ give $(1,1-\rho,-\rho)$ with $0<\rho<1$. For the Patera basis of $A_{3,8}$ one has $[u,v]=\mathcal L(u\times v)$ with $\mathcal L e_1=e_3$, $\mathcal Le_2=2e_2$, $\mathcal Le_3=e_1$, whose eigenvalues are $2,1,-1$.
    \end{proof}

    \noindent
    Since $\eta$ is left-invariant, $\tilde g=g-2\eta\otimes\eta$ is a left-invariant Lorentzian metric on the same group, and Theorem~\ref{thm:compare} applies. Table~\ref{tab:lie} lists the structures with $\psi\neq0$ for the metrics of Corollary~\ref{cor:lienormal}, up to $\xi\mapsto-\xi$ and automorphisms preserving $g$; the Ricci tensors $S$ and $\tilde S$ are computed directly from the structure constants, and in every case they satisfy $\tilde S=S+2\dv\psi\,\eta\otimes\eta$ and $\tilde r=r$, in accordance with Theorem~\ref{thm:compare}.

    \begin{table}[ht]
        \centering
        \scriptsize
        \setlength{\tabcolsep}{2.5pt}
        \caption{Left-invariant $C_{12}$-structures with $\psi\neq0$ for the metrics of Corollary~\ref{cor:lienormal}, with the Ricci tensors of $g$ and $\tilde g$; $(e^i)$ and $(f^i)$ are the dual bases of $(e_i)$ and $(f_i)$.}
        \label{tab:lie}
        \begin{tabular}{@{}lllcll@{}}
            \toprule
            Algebra & $\xi$ & $\psi$ & $\dv\psi$ & $S$ & $\tilde S$\\
            \midrule
            $A_{3,3}$ & $\cos\vartheta\,e_1+\sin\vartheta\,e_2$ & $e_3$ & $2$ & $-2g$ & $-2\tilde g$\\
            $A_{3,4}$ & $e_1$ & $e_3$ & $0$ & $-2\,e^3\otimes e^3$ & $=S$\\
            $A^\gamma_{3,5}$ & $e_1$ & $e_3$ & $\gamma+1$ & $-\mathrm{diag}(\gamma+1,\gamma(\gamma+1),\gamma^2+1)$ & $S+2(\gamma+1)\,e^1\otimes e^1$\\
            $A^\gamma_{3,5}$ & $e_2$ & $\gamma e_3$ & $\gamma(\gamma+1)$ & $-\mathrm{diag}(\gamma+1,\gamma(\gamma+1),\gamma^2+1)$ & $S+2\gamma(\gamma+1)\,e^2\otimes e^2$\\
            $A_{3,8}$ & $\frac1{\sqrt3}(e_1+e_2-e_3)$ & $-e_1-e_3$ & $0$ & $-2\,(e^1+e^3)\otimes(e^1+e^3)$ & $=S$\\
            $A_{3,8}$, $g_\rho$ & $\frac{\sqrt\rho\,f_1+f_3}{\sqrt{1+\rho}}$ & $-\sqrt\rho\,f_2$ & $0$ & $-2\rho\,f^2\otimes f^2$ & $=S$\\
            \bottomrule
        \end{tabular}
    \end{table}

    \noindent
    For $A_{3,3}$ every left-invariant metric is real hyperbolic \cite{Mil76}, and $\tilde S=-2\tilde g$ for every $\vartheta$, i.e.\ the $\eta$-reversed (Lorentzian) metric is locally anti-de Sitter space $\mathrm{AdS}_3$ (the analogue of Example~\ref{ex:S3dS}, which pairs $S^3$ with de Sitter space). In the Poincar\'e model $g=(dx^2+dy^2+dz^2)/z^2$ with $\xi=z\partial_x$, the lapse is $a=1/z$, the base is the hyperbolic plane and $\dv\psi=\Delta_ha/a=2$; Theorem~\ref{thm:compare}(ii) then gives $\tilde S=-2g+4\eta\otimes\eta=-2\tilde g$. In the unimodular cases $\dv\psi=0$, so $\tilde S=S$: the Riemannian and Lorentzian metrics have the same Ricci tensor. The entries for $A^\gamma_{3,5}$ are polynomial in $\gamma$. At the excluded value $\gamma=1$ ($A^1_{3,5}=A_{3,3}$) the two rows for $A^\gamma_{3,5}$ become the $A_{3,3}$ row with $\vartheta=0$ and $\vartheta=\pi/2$, respectively; at the excluded value $\gamma=-1$ ($A^{-1}_{3,5}=A_{3,4}$) they become the $A_{3,4}$ row, for the second one after applying the automorphism $(e_1,e_2,e_3)\mapsto(e_2,e_1,-e_3)$, which preserves $g$. The row $g_\rho$ reproduces the $A_{3,8}$ row for $\rho=\frac12$ after a homothety, with $f_2=(e_1+e_3)/\sqrt2$, and the $A_{3,4}$ row for $\rho=1$.

    \begin{proposition}[Ricci solitons]\label{prop:leftinv}
        Let $g$ be a left-invariant metric carrying a left-invariant $C_{12}$-structure.
        \begin{itemize}
            \item[(i)] The identity \eqref{eq:Tid}, with $f=\eta(V)$ constant and $W=V-f\xi$, reads $\tilde{\mathcal E}_V=\mathcal E_W+f(\eta\otimes\omega+\omega\otimes\eta)+4(\lambda+\dv\psi-\omega(W))\,\eta\otimes\eta$ for every left-invariant $V$ and every $\lambda$.
            \item[(ii)] The Ricci soliton equation for left-invariant $V$ is linear in $(V,\lambda)$. On $A_{3,3}$ its only solution, for $\tilde g$ and for $g$, is $V=0$, $\lambda=\scal_g/3$. For the structures with $\psi\neq0$ on $A_{3,4}$, $A^\gamma_{3,5}$ ($0<|\gamma|<1$) and $A_{3,8}$ the system is inconsistent, for $\tilde g$ and for $g$: there is no Ricci soliton with left-invariant potential.
            \item[(iii)] For an arbitrary smooth function $f$, $(\tilde g,f\xi,\lambda)$ is a Ricci soliton only in the following two cases: $\mathfrak g=A_{3,3}$, $S=\lambda g$ and $f\xi$ is a Killing field of $\tilde g$; or $\mathfrak g=A_{3,6}$ with a flat metric, $\psi=0$, $\lambda=0$ and $f$ constant.
        \end{itemize}
    \end{proposition}

    \begin{proof}
        (i) follows from \eqref{eq:Tid} since $df=0$. By Corollary~\ref{cor:lienormal}, in (ii) and (iii) it suffices to consider the metrics of Table~\ref{tab:lie} and the flat metric on $A_{3,6}$; automorphisms preserve left-invariance, and homotheties only rescale $\lambda$. (ii) is a computation in the bases of Table~\ref{tab:lie}; for $g_\rho$ we used the parameter $s=\sqrt\rho$, and for $g$ and $\tilde g$ the $6\times4$ coefficient matrix has rank four while the non-zero $5\times5$ minors of the augmented matrix have greatest common divisor $8s^4(s^2+1)$, resp.\ $8s^4$, so the system is inconsistent for every $\rho>0$; see also the remark below. (iii) By Theorem~\ref{thm:fxi}(i), $S=\lambda g$ on $D\times D$. By Table~\ref{tab:lie}, $S|_{D\times D}$ is not a multiple of $g|_{D\times D}$ for $A_{3,4}$, $A^\gamma_{3,5}$ (both choices of $\xi$) and $A_{3,8}$ (all $\rho$; note that $f_2\in D$ and $S$ vanishes on the unit vector of $D\cap\Span\{f_1,f_3\}$). For $A_{3,3}$ we have $S=-2g$, so $\lambda=-2=-\dv\psi$, and we are in case (B) of Theorem~\ref{thm:rigidity}. For the flat metric on $A_{3,6}$ we have $\psi=0$ and $S=0$, so $\lambda=0$ and $df=0$ by Theorem~\ref{thm:fxi}.
    \end{proof}

    \noindent
    For $A^\gamma_{3,5}$ the soliton equation is a system of six linear equations in the four unknowns $(V,\lambda)$. For $g$ and for $\tilde g$, and for both choices of $\xi$, its coefficient matrix has generic rank $4$, and the $5\times5$ minors of the augmented matrix have greatest common divisor $8\gamma(\gamma-1)(\gamma^2+1)$ up to sign, so the system could be consistent only for $\gamma\in\{0,1\}$, outside the admissible range. Ricci solitons with non-left-invariant potentials other than $f\xi$ are not excluded by Proposition~\ref{prop:leftinv}.
    \section{Corrections to the work of Delloum and Beldjilali}\label{sec:relation}
    \noindent
    Theorem~\ref{thm:compare} and its proof contradict several formulas of \cite{DB24}. Table~\ref{tab:DB} lists them together with the corrected statements; counterexamples and further details are given in Appendix~\ref{app:DB}.

    \begin{table}[ht]
        \centering
        \footnotesize
        \caption{Statements of \cite{DB24} and their corrections.}
        \label{tab:DB}
        \begin{tabular}{@{}>{\raggedright\arraybackslash}p{0.15\textwidth}>{\raggedright\arraybackslash}p{0.5\textwidth}>{\raggedright\arraybackslash}p{0.27\textwidth}@{}}
            \toprule
            In \cite{DB24} & Correct statement & Here\\
            \midrule
            (3.2) & $\nt_XY=\nabla_XY+2\eta(X)\eta(Y)\psi$ & Theorem~\ref{thm:compare}(i)\\
            (3.3) & $\tilde R(X,Y)Z=R(X,Y)Z-2\eta(Z)R(X,Y)\xi$ & Theorem~\ref{thm:compare}(v)\\
            (3.4) & $\tilde S=S+2\dv\psi\,\eta\otimes\eta$ & Theorem~\ref{thm:compare}(ii)\\
            Thm.~3.5 & $(\nt_X\varphi)Y=\eta(X)\{\omega(\varphi Y)\xi-\eta(Y)\varphi\psi\}$ & Theorem~\ref{thm:compare}(iv), Proposition~\ref{prop:lorentz}\\
            (4.1), Prop.~4.1 & $\Lie_{f\xi}\tilde g$ contains $f(\eta\otimes\omega+\omega\otimes\eta)$; a Ricci soliton $(\tilde g,\xi,\lambda)$ forces $\omega=0$ and $\lambda=0$ & Lemma~\ref{lem:Lie}(a), Corollary~\ref{cor:Vxi}\\
            Thms.~4.2, 4.4 & not correct; for $V=f\xi$ one has $df=f\omega+(\lambda+\dv\psi)\eta$, and for $V\perp\xi$ the Ricci solitons of $g$ and $\tilde g$ coincide & Theorems~\ref{thm:fxi}, \ref{thm:rigidity}, Corollary~\ref{cor:Vperp}\\
            Thm.~4.3 & not correct & Appendix~\ref{app:DB}\\
            Thm.~4.5 & not correct, both for $c_2=1$ and for $c_2=-1$ & Appendix~\ref{app:DB}\\
            p.~330 & Okumura's theorem does not apply, since $\eta\wedge d\eta=0$; with constant coefficients, $\eta$-Einstein implies Einstein or $\omega=0$ & Proposition~\ref{prop:constcoef}, Theorem~\ref{thm:nonclosed}\\
            \bottomrule
        \end{tabular}
    \end{table}

    \section{Open questions}\label{sec:open}
    \noindent
    We list the questions left open by the present work.
    \begin{itemize}
        \item[(1)] Are there compact $\eta$-Einstein para-$C_{12}$-manifolds with constant $\lambda$ and $\omega\neq0$? By Theorem~\ref{thm:compact}(a), (b), $\lambda=0$ would force $S=0$, and $\lambda\neq0$ would force $\omega$ to be exact and $e^f\xi$ to be a Killing field. In setting (R) the answer is negative by Theorem~\ref{thm:compact}(c), whose proof does not carry over to setting (P): the maximum principle and Myers' argument need the Riemannian signature, and \eqref{eq:Akappa} fails where $\nabla u$ is null.
        \item[(2)] In setting (P), describe the local structure near points of the set $W=\{d\omega\neq0\}$ at which $\Ric_h-\lambda h=0$ or all gradients $\nabla\sigma(\cdot,t)$ are null (Theorem~\ref{thm:nonclosedR}(b), Example~\ref{ex:ncP}).
        \item[(3)] Classify the $C_{12}$-manifolds with $S(X,Y)=\lambda g(X,Y)$ on $D$ for a non-constant function $\lambda$ (Remark~\ref{rem:nonconst}).
        \item[(4)] Describe the Ricci solitons of $\tilde g$ whose potential $V=f\xi+W$ has $W\neq0$ and $f\not\equiv0$ (Theorem~\ref{thm:general}), and the left-invariant para-$C_{12}$-structures on three-dimensional Lie groups with Lorentzian metrics.
    \end{itemize}
    \appendix
    \section{Details of the corrections to \texorpdfstring{\cite{DB24}}{Delloum and Beldjilali}}\label{app:DB}
    \noindent
    We give the details for Table~\ref{tab:DB}. Misprints in \cite{DB24} that do not affect the results are not listed. Following \cite{DB24}, for vector fields $V,V_1,V_2$ and constants $c_1,c_2,c_3,\lambda$ (the constant $c_3$ is denoted by $\mu$ in \cite{DB24}), the generalized $\eta$-Ricci soliton equation \cite[(2.4)]{DB24} and the generalized $\eta$-Ricci bi-soliton equation \cite[(2.5)]{DB24} read
    \begin{equation}\label{eq:gen}
    \begin{gathered}
        \Lie_Vg=-2c_1V^\flat\otimes V^\flat+2c_2S+2\lambda g+2c_3\,\eta\otimes\eta,\\
        \Lie_{V_1}g=-2c_1V_2^\flat\otimes V_2^\flat+2c_2S+2\lambda g+2c_3\,\eta\otimes\eta ,
    \end{gathered}
    \end{equation}
    and $c_3=0$ in the first equation gives the generalized Ricci soliton equation \cite[(2.2)]{DB24}.
    \begin{itemize}
        \item \cite[(3.2)]{DB24}: the Levi-Civita connections differ by $2\eta(X)\eta(Y)\psi$, not by $\eta(X)\eta(Y)\psi$ (Theorem~\ref{thm:compare}(i)). In the model, $\nt_\xi\xi-\nabla_\xi\xi=2\psi$.
        \item \cite[(3.3)]{DB24}: the correct formula is $\tilde R(X,Y)Z=R(X,Y)Z-2\eta(Z)R(X,Y)\xi$ (Theorem~\ref{thm:compare}(v)); the printed one, with coefficient $1$, fails in the model.
        \item \cite[(3.4)]{DB24}: the correct relation is $\tilde S=S+2\dv\psi\,\eta\otimes\eta$ (Theorem~\ref{thm:compare}(ii)).
        \item \cite[Thm.~3.5]{DB24}: the correct formula is $(\nt_X\varphi)Y=\eta(X)\{\omega(\varphi Y)\xi-\eta(Y)\varphi\psi\}$ (Theorem~\ref{thm:compare}(iv)); see also Proposition~\ref{prop:lorentz}.
        \item \cite[(4.1)]{DB24} and \cite[Prop.~4.1]{DB24}: the proof uses $\nt_X\xi=0$, but $\nt_X\xi=\eta(X)\psi$, and $\Lie_{f\xi}\tilde g$ contains the term $f(\eta\otimes\omega+\omega\otimes\eta)$ (Lemma~\ref{lem:Lie}(a)). By Corollary~\ref{cor:Vxi}, a Ricci soliton $(\tilde g,\xi,\lambda)$ forces $\omega=0$, $\lambda=0$ and $\tilde g$ Ricci-flat; the conclusion of \cite[Prop.~4.1]{DB24} ($\tilde g$ Einstein with $\tilde r=(2n+1)\dv\psi$) is then true only trivially ($\tilde r=0=\dv\psi$).
        \item \cite[Thms.~4.2 and 4.4]{DB24} are not correct as stated. Example~\ref{ex:S3dS} with $V=0$ (resp.\ $f\equiv0$) and $\lambda=2$ satisfies their hypotheses: $\tilde g$ is Einstein with $\tilde S=2\tilde g$, $V=0\perp\xi$, $\dv\psi=-2$ and $\lambda=\omega(V)-\dv\psi=2$. Their conclusion would be that $\Lie_\psi g=-2\omega\otimes\omega-2S+2\lambda g+2c_3\,\eta\otimes\eta$ for some $c_3$, i.e.\ \eqref{eq:gen} with $V=\psi$, $c_1=1$, $c_2=-1$. But the restriction of $\Lie_\psi g+2\omega\otimes\omega+2S-2\lambda g$ to $D\times D$ equals $-2h$, whereas $\eta\otimes\eta$ vanishes on $D$, so no $c_3$ exists. A counterexample to \cite[Thm.~4.2]{DB24} with $f\not\equiv0$ is given by Example~\ref{ex:S3dS} with $V=\partial_t$, i.e.\ $f=\cos r$: $(\tilde g,\partial_t,2)$ is a Ricci soliton, but $df=-\sin r\,dr\neq0=(\lambda+\dv\psi)\eta$, contrary to the relation $df=(\lambda+\dv\psi)\eta$ asserted there; the correct relation is $df=f\omega+(\lambda+\dv\psi)\eta$ (Theorem~\ref{thm:fxi}(ii)). The correct statements are Theorems~\ref{thm:general}, \ref{thm:fxi}, \ref{thm:rigidity} and Corollary~\ref{cor:Vperp}; for instance, for $V\perp\xi$ the Ricci solitons of $g$ and $\tilde g$ with the same potential coincide, and no bi-soliton equation is needed.
        \item \cite[Thm.~4.3]{DB24} is not correct. On $H^2\times\R$ take $g=dx^2+dy^2+e^{2x}dt^2$, i.e.\ the model \eqref{eq:model} over the flat K\"ahler plane with $\sigma=e^x$, so that $\xi=e^{-x}\partial_t$, $\psi=\partial_x$, $\omega=dx$ and $\dv\psi=|\psi|^2=1$. Here $\Lie_\psi g=2\,\eta\otimes\eta$ and $S=-dx^2-\eta\otimes\eta$, so both hypotheses of \cite[Thm.~4.3]{DB24} hold with $\lambda=0$: \eqref{eq:gen} with $V=\psi$, $c_1=1$, $c_2=-1$ and $c_3=0=2(-\lambda-\dv\psi+|\psi|^2)$, and the generalized Ricci soliton equation with $\lambda=|\psi|^2-\dv\psi=0$. However, $\tilde S=-dx^2+\eta\otimes\eta$ and $\tilde g=dx^2+dy^2-\eta\otimes\eta$, so $\tilde g$ is not $\eta$-Einstein: if $\tilde S=\alpha\tilde g+\beta\,\eta\otimes\eta$, evaluation on $\partial_y$ gives $\alpha=0$, and on $\partial_x$ gives $\alpha=-1$.
        \item \cite[Thm.~4.5]{DB24} is not correct. In Example~\ref{ex:S3dS}, $(g,0,2)$ is a Ricci soliton and $V=0\perp\xi$, but $\Lie_{-\psi}\tilde g\neq2\omega\otimes\omega+2c_2\tilde S+4\tilde g+2|\psi|^2\eta\otimes\eta$ both for $c_2=1$, as stated in \cite[Thm.~4.5]{DB24}, and for $c_2=-1$, as suggested by the computation preceding it (the two sides already differ on $D\times D$).
        \item Internal inconsistencies: in \cite[Thm.~4.2]{DB24} the value $c_3=-2(\lambda+\dv\psi-|\psi|^2)$ does not agree with the computation preceding it, which gives $c_3=-(\lambda+\dv\psi-|\psi|^2)$; in \cite[Thm.~4.5]{DB24} the value $c_2=1$ does not agree with the computation preceding it, which gives $c_2=-1$.
        \item \cite[p.~330]{DB24}: Okumura's theorem \cite{Oku62} concerns manifolds with a contact structure and cannot be applied to $C_{12}$-manifolds, which satisfy $\eta\wedge d\eta=0$. By Proposition~\ref{prop:constcoef}, $\eta$-Einstein $C_{12}$-manifolds with constant coefficients are Einstein unless $\omega=0$.
    \end{itemize}

\end{document}